\date{}
\documentclass [11pt]{article}
\usepackage{amsfonts,amsmath,graphicx}

\title{A hypergraph regularity method for generalised Tur\'an problems}
\author{
Peter Keevash \thanks{School of Mathematical Sciences,
Queen Mary, University of London, Mile End Road, London E1 4NS, UK.
Email: p.keevash@qmul.ac.uk.
Research supported in part by NSF grant DMS-0555755.}
}

\newtheorem{theo}{Theorem}[section]

\newtheorem{lemma}[theo]{Lemma}
\newtheorem{coro}[theo]{Corollary}

\newcommand{\mc}[1]{\mathcal{#1}}
\newcommand{\mb}[1]{\mathbb{#1}}
\newcommand{\nib}[1]{\noindent {\bf #1}}
\newcommand{\lra}{\leftrightarrow}
\newcommand{\sm}{\setminus}
\newcommand{\ov}{\overline}
\newcommand{\eps}{\epsilon}
\newcommand{\gen}[1]{\hspace{-1pt} \langle \hspace{-2pt} #1 \hspace{-2pt}
\rangle \hspace{-1pt}}

\newcommand{\sub}{\subseteq}

\def\COMMENT#1{}

\def\qed{\hfill $\Box$}

\begin{document}
\maketitle

\begin{abstract}
We describe a method that we believe may be foundational
for a comprehensive theory of
generalised Tur\'an problems. The cornerstone of our approach
is a quasirandom counting lemma for quasirandom hypergraphs,
which extends the standard counting lemma by not only
counting copies of a particular configuration but also showing
that these copies are evenly distributed. We demonstrate
the power of the method by proving a conjecture of Mubayi
on the codegree threshold of the Fano plane, that any
$3$-graph on $n$ vertices for which every pair of vertices
is contained in more than $n/2$ edges must contain a Fano
plane, for $n$ sufficiently large. For projective
planes over fields of odd size $q$ we show that the codegree
threshold is between $n/2-q+1$ and $n/2$, but for $PG_2(4)$
we find the somewhat surprising phenomenon that the threshold is less
than $(1/2-\eps)n$ for some small $\eps>0$.
We conclude by setting
out a program for future developments of this method
to tackle other problems.
\end{abstract}

\section{Introduction}

A famous unsolved question of Tur\'an asks for the maximum
size of a $3$-graph\footnote{A $k$-graph $G$ consists of a
vertex set $V(G)$ and an edge set
$E(G)$, each edge being some $k$-tuple of vertices.}
on $n$ vertices that does not contain a tetrahedon $K^3_4$,
i.e. $4$ vertices on which every triple is present. Despite
the superficial similarity to the analogous easy question for
graphs (the maximum size of a graph with no triangle) this
problem has evaded even an asymptotic solution for over
60 years. It may be considered a test case for the general
Tur\'an problem, that of determining the maximum size
of a $k$-graph on $n$ vertices that does not contain
some fixed $k$-graph $F$. This maximum size is called
the Tur\'an number of $F$, denoted $\mbox{ex}(n,F)$.
It is not hard to show that the limit
$\pi(F) = \lim_{n \to \infty} \mbox{ex}(n,F) / \binom{n}{k}$
exists. As a first step to understanding the Tur\'an number
we may ask to determine this limit, the Tur\'an density,
which describes the asymptotic behaviour of the Tur\'an number.
There are very few known results even in this weaker form, and
no general principles have been developed, even conjecturally.

For $0 \le s \le k$ we may define a generalised Tur\'an number
$\mbox{ex}_s(n,F)$ as the largest number $m$ such
that there is a $k$-graph $H$ on $n$ vertices that does
not contain $F$ and has minimum $s$-degree\footnote{
Given $S \sub V(H)$ the neighbourhood of $S$ in $H$
is $N_H(S) = \{T \sub V(H) \sm S: S \cup T \in E(H)\}$ and
the degree of $S$ is $|N_H(S)|$. The minimum $s$-degree
$\delta_s(H) = \min_{|S|=s} |N_H(S)|$
is the minimum of $|N_H(S)|$ over all subsets $S$ of size $s$.}
$\delta_s(H) \ge m$. Note that we recover $\mbox{ex}(n,F)$ in the
case $s=0$, and the case $s=k$ is trivial. The cases $s=0$
and $s=1$ are essentially equivalent, via a well-known
induction argument, so there is no new theory here for
graphs. However, for general hypergraphs we obtain a rich source
of new problems, and it is not apparent how they relate to
each other. There has been much recent interest
in the case $s=k-1$, which were called codegree problems
in \cite{MZ}. (See also \cite{KO,RRS} for similar questions
involving structures that are spanning rather than fixed.)
We may define generalised Tur\'an densities
as $\pi_s(F) = \lim_{n \to \infty} \mbox{ex}_s(n,F)/\binom{n-s}{k-s}$.\footnote{
In the codegree case $s=k-1$ this limit was shown to
exist in \cite{MZ}. In general the existence may be deduced from a very
general theory of Razborov \cite{R}. This is perhaps using a sledgehammer
to crack a nut, and in fact the method of \cite{MZ} can be extended
using martingale estimates in place of hypergeometric estimates. We will
elaborate slightly on the martingale aspect in the final section, but
a detailed treatment is beyond the scope we have set for this paper.}
A simple averaging argument shows that
if we define the normalised minimum $s$-degrees of $H$ as
$\hat{\delta}_s(H) = \delta_s(H)/\binom{n-s}{k-s}$ then we have
a hierarchy $\hat{\delta}_0(H) \ge \hat{\delta}_1(H)
\ge \cdots \ge \hat{\delta}_{k-1}(H)$,
so $\pi_{k-1}(F) \le \cdots \le \pi_0(F)$.

Projective geometries $PG_m(q)$ provide examples of configurations $F$ that
are surprisingly tractable for these problems. For the Fano plane
($m=q=2$) the exact Tur\'an number for $n$ sufficiently large
was determined independently and
simultaneously by Keevash and Sudakov \cite{KS1}
and F\"uredi and Simonovits \cite{FS}. They also characterised the
unique maximising configuration: a balanced complete bipartite\footnote{
A $k$-graph $H$ is bipartite if there is a partition $V(H) = A \cup B$
so that there are no edges of $H$ lying entirely within $A$ or
entirely within $B$. A complete bipartite $k$-graph contains all edges
that intersect both $A$ and $B$. It is balanced if $||A|-|B|| \le 1$.}
$3$-graph.
Earlier de Caen and F\"uredi \cite{DF}
had obtained the Tur\'an density $\pi(PG_2(2))=3/4$. On the other
hand Mubayi \cite{M1} showed that the codegree density of the Fano
plane is $\pi_2(PG_2(2)) =1/2$. He conjectured that
the exact codegree threshold satisfies $\mbox{ex}_2(n,PG_2(2)) \le n/2$.
The following result establishes this and characterises the case of equality.

\begin{theo} \label{fano}
If $n$ is sufficiently large and $H$ is a $3$-graph on $n$ vertices
with minimum $2$-degree at least $n/2$ that does not contain a Fano
plane then $n$ is even and $H$ is a balanced complete bipartite $3$-graph.
\end{theo}

General projective geometries have been studied in \cite{K} (the Tur\'an
problem) and \cite{KZ} (the codegree problem). A general bound
$\pi_q(PG_m(q)) \le 1-1/m$, was obtained in \cite{KZ},
and it was shown that equality holds whenever $m=2$ and $q$ is $2$ or odd,
and whenever $m=3$ and $q$ is $2$ or $3$. We prove the following results,
which give quite precise information about the codegree threshold for
planes over a field of odd size, and demonstrate a surprisingly different
behaviour for $PG_2(4)$.

\begin{theo} \label{pg2-odd}
Suppose $q$ is an odd prime power. Then
$\lfloor n/2 \rfloor - q + 1 \le \mbox{ex}_q(n,PG_2(q)) \le n/2$.
In the case $q=3$ and $n$ even we have $\mbox{ex}_3(n,PG_2(3)) = n/2-1$.
\end{theo}

\begin{theo} \label{pg24}
There is $\eps>0$ for which $\pi_4(PG_2(4)) < 1/2 - \eps$.
\end{theo}

The main idea in our arguments is a quasirandom counting lemma that
extends the (usual) counting lemma for quasirandom hypergraphs.
We adopt the Gowers framework as being most compatible with our
argument (there are other approaches to this theory,
see R\"odl et al. (e.g. \cite{RSk,RSc1}) and Tao \cite{Tao1}).
We will give precise definitions
later, and for now describe our result on an intuitive level.
Hypergraph regularity theory gives a method
of decomposing a hypergraph into a bounded number of pieces, each
of which behaves in an approximately random fashion. The number
of pieces depends only on the degree of approximation and is independent
of the size of the hypergraph. In order for such a decomposition
to be useful, the notion of random behaviour should be sufficiently
powerful for applications, and the general criterion that has been
used is that there should be a counting lemma, meaning a result that
a sufficiently random hypergraph contains many copies of any
small fixed configuration. Our quasirandom counting lemma will state
that not only are there many copies, but that they are uniformly
distributed within the hypergraph.

We also make use of the idea of stability,
a phenomenon which was originally discovered by Erd\H{o}s and Simonovits
in the 60's in the context of graphs with excluded subgraphs,
but has only been systematically explored relatively recently, as
researchers have realised the importance and applications of such
results in hypergraph Tur\'an theory, enumeration of discrete structures
and extremal set theory (see \cite{KM} as a recent example and for many further
references).

The rest of this paper is organised as follows.
The next section is expository in nature: it introduces the theory
needed in later sections for the special case of graphs, where it will
be mostly familiar to many readers (although our quasirandom counting lemma
is new even for graphs). Then in section 3 we introduce the Gowers quasirandomness
framework for $3$-graphs and present a case of our quasirandom counting lemma
that we will need to prove Theorem \ref{fano}. Section 4 contains the
proof of Theorem \ref{fano}, using the quasirandomness theory from section 3
and also the method of `stability', or approximate structure.
In section 5 we present the general theory of quasirandomness hypergraphs and the
full form of our quasirandom counting lemma: this is the engine behind our entire
approach. We also give an application to
generalised Tur\'an problem for configurations that have a certain
special form. This general theory is applied in section 6 to the
study of codegree problems in projective planes, where we prove the other theorems
stated above. The final section sets out a program
for future developments of this method to other generalised Tur\'an problems.
Since our formulation of the Gowers quasirandomness framework uses some
non-trivial variations on the original framework,
we give justifications for these variations in an appendix to the paper.

\nib{Notation.}\
Write $[n]=\{1,\cdots,n\}$. If $X$ is a set and $k$ is a number
then $\binom{X}{k} = \{Y \sub X: |Y|=k\}$,
$\binom{X}{\le k} = \cup_{i \le k} \binom{X}{i}$ and
$\binom{X}{< k} = \cup_{i < k} \binom{X}{i}$.
$a \pm b$ denotes an unspecified real number in the interval
$[a-b,a+b]$. It is convenient to
regard a finite set $X$ as being equipped with
the uniform probability measure $\mb{P}(\{x\})=1/|X|$,
so that we can express the average of a function $f$ defined on
$X$ as $\mb{E}_{x \in X} f(x)$. A $k$-graph $H$ consists of a
vertex set $V(H)$ and an edge set
$E(H)$, each edge being some $k$-tuple of vertices.
We often identify $H$ with $E(H)$, thus $|H|$ is the number
of edges in $H$.
Given $S \sub V(H)$ the neighbourhood of $S$ in $H$
is $N_H(S) = \{T \sub V(H) \sm S: S \cup T \in E(H)\}$ and
the degree of $S$ is $|N_H(S)|$. The minimum $s$-degree
$\delta_s(H) = \min_{|S|=s} |N_H(S)|$
is the minimum of $|N_H(S)|$ over all subsets $S$ of size $s$.
Given $X \sub V(H)$ the restriction $H[X]$ is a $k$-graph
with vertex set $X$ and edge set equal to all those edges
of $H$ that are contained in $X$.
Suppose $F$ and $H$ are $k$-graphs. The homomorphism density
$d_F(H)$ is the probability that a randomly chosen map
$\phi:V(F) \to V(H)$ is a homomorphism, i.e. $\phi(e)$ is an
edge of $H$ for every edge $e$ of $F$ (which we
also write as $\phi(F) \sub H$).
We also use the same notation in a `partite setting' (this will be
explained when it occurs).
When $F=e$ consists of
just a single edge we write $d(H)=d_e(H)=k!|E(H)||V(H)|^{-k}$,
and call this the density of $H$.
We use the notation $0 < \alpha \ll \beta$ to mean that there is an increasing function
$f(x)$ so that the ensuing argument is valid for $0 < \alpha < f(\beta)$.

\section{Graphs: regularity and counting, quasirandomness and quasirandom counting.}

The purpose of this section is expository: we introduce the theory needed
in later sections for the special case of graphs, where it is considerably
simpler, and partly familiar to many readers.

In the first subsection we describe Szemer\'edi's regularity lemma \cite{Sz}, one of the most powerful tools in
modern graph theory. Roughly speaking, it says that any graph can be approximated by an average
with respect to a partition of its vertex set into a bounded number of
classes, the number of classes depending only on the accuracy of the
desired approximation, and not on the number of vertices in the
graph. Each pair of classes span a bipartite subgraph that is `regular',
meaning that the proportion of edges in any large bipartite subgraph
is close to the proportion of edges in the pair as a whole.
A key property of this approximation is that it leads to a
`counting lemma', allowing an accurate prediction of the number of
copies of any small fixed graph spanned by some specified classes of the
partition. We refer the reader to \cite{KoS} for a survey of the
regularity lemma and its applications.

The second subsection discusses quasirandomness of graphs, a concept
introduced by Chung, Graham and Wilson \cite{CGW} (see also Thomason \cite{Th}
for a similar notion). There are many ways of describing this concept, all of
which are broadly equivalent (up to renaming constants); in fact, it is
also equivalent to regularity (as described in the first subsection).
A particularly simple formulation is to call a bipartite graph quasirandom
if the number of $4$-cycles is close to what would be expected in a random
graph with the same edge density. A closely related formulation that
forms the basis for the Gowers approach to quasirandomness in hypergraphs
is to say that if we count $4$-cycles weighted by the `balanced function' of
the graph then the result is small. Our discussion in this subsection
is based on section 3 of \cite{G1} (we are more brief on those points discussed
there, but we also provide some additional arguments that are omitted there).

In the third subsection we introduce the graph case of our quasirandom counting lemma,
an extension of the counting lemma discussed in the first subsection,
saying that copies of any small fixed graph are well-distributed in
the graph. This is a new result even in the special case of graphs,
and has consequences that are somewhat surprising at first sight.

\subsection{Regularity and counting}

We start by describing the notion of `regularity' for bipartite graphs.
The density of a bipartite graph $G = (A,B)$
with vertex classes $A$ and $B$ is defined to be
\[d_G(A,B) := \frac{e_G(A,B)}{|A||B|}.\]
We often write $d(A,B)$ if this is unambiguous. Given $\eps>0$, we say
that $G$ is \emph{$\eps$-regular} if
for all subsets $X\subseteq A$ and $Y\subseteq B$
with $|X|>\eps|A|$ and $|Y|>\eps|B|$ we have that
$|d(X,Y)-d(A,B)|<\eps$.

The regularity lemma says that any graph can be partitioned into a bounded
number of regular pairs and a few leftover edges. Formally:

\begin{theo}
For every real $\eps>0$ and number $m_0 \ge 1$ there are numbers
$m, n_0 \ge 1$ so that for any graph $G$ on $n \ge n_0$ vertices
we can partition its vertices as $V(G) = V_0 \cup V_1 \cup \cdots \cup V_k$
so that
\begin{itemize}
\item $m_0 \le k \le m$,
\item $|V_0| < \eps n$,
\item $|V_1|=|V_2|=\cdots=|V_k|$, and
\item $(V_i,V_j)$ spans an $\eps$-regular bipartite subgraph of $G$
for all but at most $\eps k^2$ pairs $1 \le i < j \le k$.
\end{itemize}
\end{theo}

\nib{Remarks.}
In applications one takes $k \ge m_0 \ge \eps^{-1}$ so that the number
of edges within any $V_i$ is negligible.
An `exceptional class' $V_0$ is allowed so that the remaining
partition can be `equitable', i.e. $|V_1|=|V_2|=\cdots=|V_k|$. If one
prefers not to have an exceptional class then its vertices may be distributed
among the other classes to obtain a partition with the same regularity
properties for a slightly larger $\eps$ with the class sizes
differing by at most one. We refer the reader to Section 1 of
\cite{KoS} for further discussion of variants of the regularity lemma.

\medskip

A key property of regularity is that we can accurately count copies of
any small fixed graph. For the purpose of exposition we state two simple
cases, and then the general `counting lemma'. Throughout we assume a hierarchy
$0 < 1/n \ll \eps \ll d$, that every density we consider is at least $d$
and each part in our graphs contains at least $n$ vertices. The following
well-known statements can be proved using similar arguments to that given
for Lemma 2.1 in \cite{KoS}.
\begin{description}
\item[Triangles.]
Suppose $G$ is a tripartite graph with parts
$V_1$, $V_2$, $V_3$ and each pair $(V_i,V_j)$ spans an $\eps$-regular
bipartite graph of density $d_{ij}$. Let $\triangle(G)$ be the set of triangles in $G$.
Then we can estimate the `triangle density' in $G$ as
$d_{\triangle}(G) = |\triangle(G)|/|V_1||V_2||V_3| = d_{12}d_{13}d_{23} \pm 8\eps$.
\footnote{The constant $8$ is not best possible, but we only care that the error
should tend to zero as $\eps$ tends to $0$.}
\item[$4$-cycles.]
Suppose $G=(X,Y)$ is an $\eps$-regular bipartite
graph with density $d$. Let $C_4(G)$ be the number of labelled $4$-cycles in $G$,
i.e. quadruples $(x_1,x_2,y_1,y_2)$ with $x_1 \ne x_2 \in X$, $y_1 \ne y_2 \in Y$
such that $x_1y_1$, $x_1y_2$, $x_2y_1$ and $x_2y_2$ are all edges of $G$.
We may define $d_{C_4}(G)$, the `bipartite homomorphism density' of $C_4$ in $G$,
as follows. Fix a $4$-cycle $C_4=(A,B)$, considered as a bipartite graph $K_{2,2}$
on $A=\{a_1,a_2\}$ and $B=\{b_1,b_2\}$. Let $\Phi$ be the set of all
`bipartite maps' $\phi$ from $A \cup B$ to $X \cup Y$,
i.e. functions with $\phi(A) \sub X$ and $\phi(B) \sub Y$.
Define
\footnote{By $\phi(C_4) \sub G$ we mean $\phi(e) \in E(G)$ for each edge $e$ of $C_4$.}
$$d_{C_4}(G) = \mb{P}_{\phi \in \Phi}[\phi(C_4) \sub G]
= C_4(G)/|X|^2|Y|^2 \pm O(1/n).$$
Then $d_{C_4}(G) = d^4 \pm 10\eps$.
In fact, a lower bound of $d^4$ follows from the Cauchy-Schwartz inequality,
so $\eps$-regularity of $G$ is only needed to prove the upper bound.
\item[General graphs.]
Suppose $G$ is an $r$-partite graph on $V = V_1 \cup \cdots \cup V_r$ and
each pair $(V_i,V_j)$, $1 \le i<j \le r$ spans an $\eps$-regular bipartite
graph with density $d_{ij}(G)$. Suppose $H$ is an $r$-partite graph
on $Y = Y_1 \cup \cdots \cup Y_r$. Let $\Phi(Y,V)$ be the set of
$r$-partite maps from $Y$ to $V$, i.e. maps $\phi:Y \to V$ with
$\phi(Y_i) \sub V_i$ for $1 \le i \le r$.
Define the $r$-partite homomorphism density of
$H$ in $G$ as $d_H(G) = \mb{P}_{\phi \in \Phi(Y,V)}[\phi(H) \sub G]$.
Then $d_H(G) = \prod_{e \in E(H)} d_e(G) \pm O_H(\eps)$,
where $d_e(G)$ means that density $d_{ij}(G)$ for which $e \in H[Y_i,Y_j]$
and the $O_H$-notation indicates that the implied constant depends only on $H$.
\end{description}

\subsection{Quasirandomness}

The regularity property discussed in the previous subsection turns out
to be characterised by the counting lemma; in fact, somewhat surprisingly,
it is characterised just by counting $4$-cycles. To be precise, if
$0 < 1/n \ll \eps \ll \eps' \ll d \le 1$ and $G=(X,Y)$ is a bipartite
graph with $|X|, |Y| \ge n$, density $d$ and $C_4$-density $d_{C_4}(G) < d^4 + \eps$
then $G$ is $\eps'$-regular. In order to illuminate some later more general
arguments we will prove this fact here,
via a closely related characterisation of Gowers (see section 3
of \cite{G1}) that forms the basis for the Gowers approach to quasirandomness
in hypergraphs.

Suppose $G$ is a bipartite graph with parts $X$ and $Y$. We can identify $G$
with its characteristic function $G(x,y)$, which for $x \in X$ and $y \in Y$
is defined to be $1$ if $xy$ is an edge of $G$, or $0$ otherwise. Define
the {\em balanced function} $\ov{G}(x,y) = G(x,y)-d$, where
$d=d_G(X,Y)=|E_G(X,Y)|/|X||Y|$ is the density of $G$,
i.e. $\ov{G}(x,y)$ is $1-d$ if $xy$ is an edge of $G$, or $-d$ otherwise.
Note that $\sum_{x \in X, y \in Y} \ov{G}(x,y)=0$. For any function
$f$ defined on $X \times Y$ we define
$$C_4(f) = |X|^{-2}|Y|^{-2}\sum_{x_1,x_2 \in X}\sum_{y_1,y_2 \in Y}
f(x_1,y_1)f(x_1,y_2)f(x_2,y_1)f(x_2,y_2).$$
We can rephrase this by recalling the setup in the previous subsection,
where we had a $4$-cycle $C_4=(A,B)$, considered as a bipartite graph $K_{2,2}$
on $A=\{a_1,a_2\}$ and $B=\{b_1,b_2\}$, and let $\Phi$ denote the set of all
bipartite maps $\phi$ from $A \cup B$ to $X \cup Y$. Then
$$C_4(f) = \mb{E}_{\phi \in \Phi} \prod_{e \in E(C_4)} f(\phi(e)).$$
In particular $C_4(G)=d_{C_4}(G)$, where the first instance of $G$ is to
be understood as the characteristic function $G(x,y)$.

We say that a function $f: X \times Y \to [-1,1]$ is {\em $\eta$-quasirandom} if
$C_4(f) < \eta$, and we say that $G$ is {\em $\eta$-quasirandom} if its balanced
function $\ov{G}$ is $\eta$-quasirandom.

\begin{theo} \label{c4-tfae}
Suppose $G$ is a bipartite graph of density $d$
with parts $X$ and $Y$ of size at least $n$.
The following are equivalent, in the sense that each implication $(i)\Rightarrow(j)$
is true for $0 < 1/n \ll \eps_i \ll \eps_j \ll d \le 1$:

(1) $G$ is $\eps_1$-regular.

(2) $d_{C_4}(G) = d^4 \pm \eps_2$.

(3) $C_4(\ov{G}) < \eps_3$.

\end{theo}

Before giving the proof we quote a simple version of the `second moment method',
Lemma 6.5 in \cite{G1}: if $0< \alpha,d < 1$,
$\sum_{i=1}^n a_i \ge (d-\alpha)n$ and
$\sum_{i=1}^n a_i^2 \le (d^2+\alpha)n$ then $a_i = d \pm \alpha^{1/4}$ for
all but at most $3\alpha^{1/2}n$ values of $i$. Indeed,
$\sum_{i=1}^n (a_i-d)^2 = \sum a_i^2 - 2d \sum a_i + d^2 n
\le (d^2 + \alpha - 2d(d-\alpha) + d^2)n < 3\alpha n$.
Also, we often use the Cauchy-Schwartz inequality in the form
$0 \le \mb{E}(Z - \mb{E}Z)^2 = \mb{E}Z^2 - (\mb{E}Z)^2$, for a random variable $Z$.

\medskip

\nib{Proof.}
$(1)\Rightarrow(2)$: This implication is given by the counting lemma quoted in the
previous subsection, but for completeness we give a proof here. Consider a
random map $\phi \in \Phi$, i.e. a random bipartite map from $A \cup B$ to $X \cup Y$,
where as before we consider $C_4$ as a bipartite graph with parts
$A=\{a_1,a_2\}$ and $B=\{b_1,b_2\}$. Let $E_1$ be the event that
$x_1 = \phi(a_1)$ has $(d \pm \eps_1)|Y|$ neighbours in $Y$. By definition
of $\eps_1$-regularity we have $\mb{P}(E_1) > 1 - 2\eps_1$ (there are
at most $\eps_1 |X|$ vertices with more than $(d + \eps_1)|Y|$ neighbours and
at most $\eps_1 |X|$ vertices with less than $(d - \eps_1)|Y|$ neighbours).
Let $E_2$ be the event that $x_2 = \phi(a_2)$ has $(d \pm \eps_1)^2|Y|$ neighbours in
$N(x_1)$. Again, $\eps_1$-regularity gives $\mb{P}(E_2|E_1) > 1 - 2\eps_1$.
Now the event $\{\phi(C_4) \sub G\}$ occurs if and only if $y_1 = \phi(b_1)$ and
$y_2 = \phi(b_2)$ lie in $N(x_1) \cap N(x_2)$, so we have
$$d_{C_4}(G) = \mb{P}_{\phi \in \Phi}[\phi(C_4) \sub G]
= \mb{P}_{\phi \in \Phi}[\{\phi(C_4) \sub G\} \cap E_1 \cap E_2] \pm 4\eps_1
= d^4 \pm 10\eps_1.$$
This proves the implication with $\eps_2 = 10\eps_1$.

$(2)\Rightarrow(3)$: Suppose that $d_{C_4}(G) = d^4 \pm \eps_2$.
Our first step is to show that we can also count any subgraph of $C_4$,
in that $d_{P_2}(G) = d^2 \pm \eps_2$ and $d_{P_3}(G) = d^3 \pm 4\eps_2^{1/4}$,
where $P_i$ is the path with $i$ edges (we already know that
$d_{P_1}(G)=d_G(X,Y)=d$, and it is immediate that for a matching
$M_2$ of two edges we have $d_{M_2}(G)=d^2$).
We start with the homomorphism density of $P_2$,
say with the central vertex being mapped to $X$ and the two outer vertices
to $Y$ (the same bound will hold vice versa). We can write
$$d_{P_2}(G) = \mb{E}_{x \in X, y_1 \in Y, y_2 \in Y} G(x,y_1)G(x,y_2)
= \mb{E}_{x \in X} (\mb{E}_{y \in Y} G(x,y))^2,$$
which by Cauchy-Schwartz is at least
$(\mb{E}_{x \in X,y \in Y} G(x,y))^2=d^2$.
On the other hand, we can again apply Cauchy-Schwartz to get
\begin{align*}
d_{P_2}(G)^2 & = (\mb{E}_{x, y_1, y_2} G(x,y_1)G(x,y_2))^2
\le \mb{E}_{y_1,y_2} (\mb{E}_x G(x,y_1)G(x,y_2))^2 \\
& = \mb{E}_{y_1,y_2} \mb{E}_{x_1,x_2} G(x_1,y_1)G(x_1,y_2)G(x_2,y_1)G(x_2,y_2)
= d_{C_4}(G) < d^4 + \eps_2.
\end{align*}
This gives $d_{P_2}(G) = d^2 \pm \eps_2$, and moreover,
applying the second moment method quoted before the proof,
for a random $x \in X$,
with probability at least $1 - 3\eps_2^{1/2}$ we have
$\mb{E}_{y_1,y_2} G(x,y_1)G(x,y_2) = d^2 \pm \eps_2^{1/4}$,
so $|N(x)|/|Y| = \mb{E}_y G(x,y) = d \pm \eps_2^{1/4}$.
This allows us to estimate
$$d_{P_3}(G) = \mb{E}_{x_1,x_2 \in X, y_1 \in Y, y_2 \in Y} G(x_1,y_1)G(x_1,y_2)G(x_2,y_2)
= (d^2 \pm \eps_2^{1/4})(d \pm \eps_2^{1/4}) \pm 6\eps_2^{1/2}
= d^3 \pm 4\eps_2^{1/4},$$
where the main term gives the contribution when $x_1$ and $y_2$ have typical
neighbourhoods and $6\eps_2^{1/2}$ bounds the error coming from atypical $x_1$ and $y_2$.

Now we can estimate $C_4(\ov{G})$. Write $\ov{G} = f_0 - f_1$, where
$f_0(x,y)=G(x,y)$ and $f_1(x,y)=d$ (a constant function). Then
\begin{align*}
C_4(\ov{G}) & = \mb{E}_{x_1,x_2 \in X, y_1,y_2 \in Y}
\ov{G}(x_1,y_1)\ov{G}(x_1,y_2)\ov{G}(x_2,y_1)\ov{G}(x_2,y_2) \\
& = \sum_M (-1)^{\sum M} \mb{E}_{x_1,x_2,y_1,y_2}
f_{M_{11}}(x_1,y_1)f_{M_{12}}(x_1,y_2)f_{M_{21}}(x_2,y_1)f_{M_{22}}(x_2,y_2),
\end{align*}
where $M=(M_{ij})_{i,j \in \{1,2\}}$ ranges over $2 \times 2$ matrices
with $\{0,1\}$-entries and $\sum M$ is $\sum_{i,j} M_{ij}$.
For each $M$ we can estimate the summand corresponding to $M$,
using the estimate for the homomorphism density of the subgraph of
$C_4$ corresponding to the $0$-entries of $M$. For example,
if $M_{11}=M_{12}=M_{21}=0$ and $M_{22}=1$, then $M$ corresponds to a $P_3$,
and the corresponding summand is $- d \cdot d_{P_3}(G) = - d^4 \pm 4\eps_2^{1/4}$.
Thus we obtain $8$ summands in the range $d^4 \pm 4\eps_2^{1/4}$
and $8$ summands in the range $- d^4 \pm 4\eps_2^{1/4}$,
giving a total of at most $64\eps_2^{1/4}$. This proves
the implication with $\eps_3 = 64\eps_2^{1/4}$.

$(3)\Rightarrow(1)$: Suppose $C_4(\ov{G}) < \eps_3$. Consider
$X' \sub X$ with $|X'| > \eps_1 |X|$ and $Y' \sub Y$ with $|Y'| > \eps_1 |Y|$.
Write $X'(x)$ for the characteristic function of $X'$, i.e. $X'(x)$ is $1$
if $x \in X'$, otherwise it is $0$; similarly, let $Y'(y)$ be the characteristic
function of $Y'$. Then
$$d_G(X',Y')-d = \mb{E}_{x \in X', y \in Y'} \ov{G}(x,y)
= \frac{|X||Y|}{|X'||Y'|} \mb{E}_{x \in X, y \in Y} X'(x)Y'(y)\ov{G}(x,y)$$
and by Cauchy-Schwartz (in the first and third inequalities below) we have
\begin{align*}
(\mb{E}_{x \in X, y \in Y} X'(x)Y'(y)\ov{G}(x,y))^4
& \le (\mb{E}_{x \in X} (\mb{E}_{y \in Y} X'(x)Y'(y)\ov{G}(x,y))^2)^2 \\
& \le (\mb{E}_{x \in X} (\mb{E}_{y \in Y} Y'(y)\ov{G}(x,y))^2)^2 \\
& = (\mb{E}_{x \in X} \mb{E}_{y_1,y_2 \in Y} Y'(y_1)Y'(y_2)\ov{G}(x,y_1)\ov{G}(x,y_2))^2 \\
& \le \mb{E}_{y_1,y_2 \in Y} (\mb{E}_{x \in X} Y'(y_1)Y'(y_2)\ov{G}(x,y_1)\ov{G}(x,y_2))^2 \\
& \le \mb{E}_{y_1,y_2 \in Y} (\mb{E}_{x \in X} \ov{G}(x,y_1)\ov{G}(x,y_2))^2 \\
& =  \mb{E}_{y_1,y_2 \in Y} \mb{E}_{x_1,x_2 \in X} \ov{G}(x_1,y_1)\ov{G}(x_1,y_2)
\ov{G}(x_2,y_1)\ov{G}(x_2,y_2) \\
& = C_4(\ov{G}) < \eps_3.
\end{align*}
This gives $d_G(X',Y') = d \pm \eps_1^{-2} \eps_3^{1/4}$, which proves the
implication with $\eps_1 = \eps_3^{1/12}$. \qed

\subsection{Quasirandom counting}

Now we will introduce the graph case of our quasirandom counting lemma,
an extension of the counting lemma discussed in the first subsection,
saying that copies of any small fixed graph are well-distributed in
the graph. As before, for the purpose of exposition we lead up to the
general case through two illustrative cases.
\begin{description}
\item[Triangles.]
Suppose $G$ is a tripartite graph with parts
$V_1$, $V_2$, $V_3$ and each pair $(V_i,V_j)$ spans an $\eps$-regular
bipartite graph of density $d_{ij}$.
We remarked before that the `triangle density' in $G$ can be estimated as
$d_{\triangle}(G) = d_{12}d_{13}d_{23} \pm 8\eps$.
Moreover, it is easy to see that the triangles of $G$ are `well-distributed', in that
any sufficiently large subsets $V'_i \sub V_i$, $1 \le i \le 3$ induce a subgraph of
$G$ that also has triangle density about $d_{12}d_{13}d_{23}$. Indeed, if
$V'_i \sub V_i$ with $|V'_i| > \eps^{1/2}|V_i|$, $1 \le i \le 3$, then each
$(V'_i,V'_j)$ induces an $\eps^{1/2}$-regular subgraph of $G$ with
density $d_{ij} \pm \eps$, so we can apply triangle counting directly to
see that $(V'_1,V'_2,V'_3)$ has triangle density
$(d_{12} \pm \eps)(d_{13} \pm \eps)(d_{23} \pm \eps) \pm 8\eps^{1/2}
= d_{12}d_{13}d_{23} \pm 15\eps^{1/2}$ (say).
\item[$4$-cycles.]
Suppose $G=(X,Y)$ is an $\eps$-regular bipartite graph with density $d$.
We saw before that we can estimate the $C_4$-density of $G$ as $d_{C_4}(G) = d^4 \pm 10\eps$.
A similar argument to that just given for triangles allows us to estimate
the $C_4$-density of $G$ restricted to sufficient large subsets $X_0 \sub X$
and $Y_0 \sub Y$. But in fact, we can make a stronger claim: given any sufficiently
dense graphs $H_X$ on $X$ and $H_Y$ on $Y$ we can estimate the
density of $4$-cycles $\{x_1y_1,x_1y_2,x_2y_1,x_2y_2\}$ of $G$ in
which $x_1x_2$ is an edge of $H_X$ and $y_1y_2$ is an edge of $H_Y$.

To do this, first recall our definition of
$d_{C_4}(G)$ as $\mb{P}_{\phi \in \Phi}[\phi(C_4) \sub G]$,
where we fix a $4$-cycle $C_4=(A,B)$, considered as a bipartite graph
on $A=\{a_1,a_2\}$ and $B=\{b_1,b_2\}$ and let $\Phi$ be the set of all
bipartite maps $\phi$ from $A \cup B$ to $X \cup Y$.
Consider an auxiliary `$C_4$-homomorphism' graph $G'$, defined as follows.
$G'$ is a bipartite graph with parts $X'$ and $Y'$, where
$X'$ consists of all maps $\phi_1:A \to X$ and
$Y'$ consists of all maps $\phi_2:B \to Y$.
Given $\phi_1 \in X'$ and $\phi_2 \in Y'$ we can construct
a bipartite map $\phi = (\phi_1,\phi_2) \in \Phi$ in an obvious
manner: $\phi$ restricts to $\phi_1$ on $A$ and $\phi_2$ on $B$.
We say that $\phi_1\phi_2$ is an edge of $G'$ if
$\phi=(\phi_1,\phi_2)$ is a homomorphism from $C_4$ to $G$,
i.e. $\phi(a_i,b_j) \in E(G)$ for $i,j \in \{1,2\}$.

We claim that if $0 < \eps \ll \eps'$ and $G$ is $\eps$-regular then $G'$ is $\eps'$-regular.
To see this we use Theorem \ref{c4-tfae}. By construction $G'$ has density
$d_{G'}(X',Y')=d_{C_4}(G)=d^4 \pm 10\eps$. Also we can calculate
$d_{C_4}(G')=d_{K_{4,4}}(G)$, where $K_{4,4}$ is the complete bipartite
graph with $4$ vertices in each part. To see this, label the parts of
the $K_{4,4}$ as $A^1 \cup A^2$ and $B^1 \cup B^2$,
where $A^1 = \{a_1^1,a_2^1\}$, $A^2 = \{a_1^2,a_2^2\}$,
$B^1 = \{b_1^1,b_2^1\}$, $B^2 = \{b_1^2,b_2^2\}$.
Then we can identify a bipartite map $\psi'$ from $A \cup B$ to $X' \cup Y'$
with a bipartite map $\psi$ from $(A^1 \cup A^2) \cup (B^1 \cup B^2)$ to $X \cup Y$:
define $\psi(a_i^j) = \psi'(a_i)(a_j)$ and $\psi(b_i^j) = \psi'(b_i)(b_j)$.
\footnote{To interpret this notation, observe that e.g. $\psi'(a_1) = \omega \in X'$
is a map from $A$ to $X$, so $\psi(a_1^2)=\omega(a_2) \in X$.}
Now $\psi'$ is a homomorphism from $C_4$ to $G'$ if and only if
$\psi$ is a homomorphism from $K_{4,4}$ to $G$:
we have $\{\psi'(a_i),\psi'(b_k)\} \in G'$ $\forall i,k \in \{1,2\}$
$\lra$ $\{\psi'(a_i)(a_j),\psi'(b_k)(b_\ell)\} \in G$ $\forall i,j,k,\ell \in \{1,2\}$
$\lra$ $\{\psi(a_i^j),\psi(b_k^\ell)\} \in G$ $\forall i,j,k,\ell \in \{1,2\}$.
Therefore $d_{C_4}(G')=d_{K_{4,4}}(G)$.
But we can estimate $d_{K_{4,4}}(G) = d^{16} \pm O(\eps)$ by the counting lemma,
so $d_{C_4}(G')=d_{G'}(X',Y')^4 \pm O(\eps)$. Now Theorem \ref{c4-tfae}
tells us that $G'$ is $\eps'$-regular.

Finally, suppose that we have graphs $H_X$ on $X$ and $H_Y$ on $Y$
with $|E(H_X)| > \eps'|X|^2$ and $|E(H_Y)| > \eps'|Y|^2$. Applying the
definition of $\eps'$-regularity, we see that if we choose a random
edge $x_1x_2$ of $H_X$ and (independently) a random edge $y_1y_2$ of $H_Y$ then
$\{x_1y_1,x_1y_2,x_2y_1,x_2y_2\}$ is a $4$-cycle in $G$ with
probability $d^4 \pm 2\eps'$.

\item[General graphs.]
Suppose $G$ is an $r$-partite graph on $V = V_1 \cup \cdots \cup V_r$ and
each pair $(V_i,V_j)$, $1 \le i<j \le r$ spans an $\eps$-regular bipartite
graph with density $d_{ij}(G)$. Suppose $H$ is an $r$-partite graph
on $Y = Y_1 \cup \cdots \cup Y_r$. Let $\Phi(Y,V)$ be the set of
$r$-partite maps from $Y$ to $V$, i.e. maps $\phi:Y \to V$ with
$\phi(Y_i) \sub V_i$ for $1 \le i \le r$.
Recall that the $r$-partite homomorphism density of
$H$ in $G$ is $d_H(G) = \mb{P}_{\phi \in \Phi(Y,V)}[\phi(H) \sub G]$.
The counting lemma says that
$d_H(G) = \prod_{e \in E(H)} d_e(G) \pm O_H(\eps)$,
where $d_e(G)$ means that density $d_{ij}(G)$ for which $e \in H[Y_i,Y_j]$
and the $O_H$-notation indicates that the implied constant depends only on $H$.

Moreover, copies of $H$ are well-distributed in $G$ in the following sense.
Consider the auxiliary `$H$-homomorphism' graph $G'$, defined as follows.
$G'$ is an $r$-partite graph on $V' = V'_1 \cup \cdots \cup V'_r$,
where $V'_i$ consists of all maps $\phi_i:Y_i \to V_i$.
Given $\phi_i \in V'_i$ and $\phi_j \in V'_j$, for some $1 \le i < j \le r$
we say that $\phi_i\phi_j$ is an edge of $G'$ if $\phi_i(y_i)\phi_j(y_j)$
is an edge of $G$ for every edge $y_iy_j$ of $H$ with $y_i \in Y_i$
and $y_j \in Y_j$. Then copies of the complete graph $K_r$ in $G'$
correspond to homomorphisms from $H$ to $G$.

Suppose we have a constant hierarchy $0 < \eps \ll \eps' \ll \eps'' < 1$. A similar argument
to that used in the case of $4$-cycles shows that for every $1 \le i < j \le r$,
$(V'_i,V'_j)$ induces a subgraph of $G'$ with density
$d_{ij}^{e_H(Y_i,Y_j)} + O(\eps)$ that is $\eps'$-regular.
Suppose also that we have a $|Y_i|$-graph $J_i$ on $V_i$
with $|E(J_i)| > \eps'' |V_i|^{|Y_i|}$ for $1 \le i \le r$.
Then by $\eps'$-regularity of $G'$ and the counting lemma for $K_r$,
we see that if choose independent random edges $e_i$ in $J_i$ for $1 \le i \le r$
then $\{e_1,\cdots,e_r\}$ spans a copy of $H$ in $G$ with probability
$\prod_{e \in E(H)} d_e(G) \pm \eps''$.
\end{description}

\section{Quasirandom $3$-graphs}

In this section we discuss the Gowers approach to quasirandomness in $3$-graphs:
our exposition will be quite condensed, and for more details we
refer the reader to that given in \cite{G1}. Although we will later
repeat this discussion for general $k$-graphs, we feel it is helpful to first present
the case $k=3$, which is simpler to grasp for a reader new to the subject.
We conclude this section by describing our quasirandom counting lemma for the
special case of the distribution of octahedra in $3$-graphs, which is analogous
to the distribution of $4$-cycles in graphs described above;
our solution of Mubayi's conjecture, Theorem \ref{fano}, will only make use of this case.

\subsection{Two cautionary examples}

When considering how to generalise regularity from graphs to $3$-graphs,
a natural first attempt is to take a $3$-graph $H$ and partition its vertices as
$V(H) = V_1 \cup \cdots \cup V_k$, for some $k(\eps)$, so that all but
at most $\eps k^3$ triples $(V_a,V_b,V_c)$ span a tripartite
$3$-graph that is {\em $\eps$-vertex-regular}, meaning that
for any $V'_a \sub V_a$, $V'_b \sub V_b$, $V'_c \sub V_c$
with $|V'_a| \ge \eps|V_a|$, $|V'_b| \ge \eps|V_b|$, $|V'_c| \ge \eps|V_c|$
we have $\mb{E}_{x \in V'_a,y \in V'_b,z \in V'_c} H(x,y,z)
= \mb{E}_{x \in V_a,y \in V_b,z \in V_c} H(x,y,z) \pm \eps$.
This is indeed possible, as shown by Chung and Graham \cite{CG},
with a proof closely modelled on that of the graph regularity lemma.
This result is often known as the `weak hypergraph regularity lemma',
as although it does have some applications, the property of vertex-regularity
is not strong enough to prove a counting lemma, as the following example
of R\"odl demonstrates quite dramatically.

\medskip

\nib{Example.} Take a set $X = X_1 \cup X_2 \cup X_3 \cup X_4$,
where $X_i$, $1 \le i \le 4$ are pairwise disjoint sets of size $n$.
Consider a random orientation of the complete $4$-partite graph on $X$,
i.e. for every $x_i \in X_i$, $x_j \in X_j$, $1 \le i < j \le 4$
we choose the arc $x_ix_j$ or the arc $x_jx_i$, each choice having
probability $1/2$, all choices being independent.
Define a $3$-graph $H$ on $X$ to consist of all triples $x_ix_jx_k$
that induce a cyclic triangle in the orientation. Then with high
probability each triple $(X_i,X_j,X_k)$ spans an $\eps$-vertex-regular
triple when $n \gg \eps^{-1}$, and this will remain true even if we partition
$H$ into $k(\eps)$ parts. However, it is easy to see that $H$ does not
contain any copy of the tetrahedron $K_4^3$.

\medskip

In the light of this example, we might informally say that vertex-regularity
uses a random model of a $3$-graph in which triples are randomly chosen as edges
with some probability $p$, but that a counting lemma fails because of potential
correlations between edges that share a pair of vertices. Our second example
(adapted from \cite{G1}) points the way to a better model.

\medskip

\nib{Example.} Take a set $X = X_1 \cup X_2 \cup X_3$,
where $X_i$, $1 \le i \le 3$ are pairwise disjoint sets of size $n$.
For each pair $1 \le i < j \le 3$ define a random bipartite graph
$G_{ij}$, in which each $x_ix_j$, $x_i \in X_i$, $x_j \in X_j$ is
an edge with probability $d_{ij}$, all choices being independent.
Consider $\triangle = \triangle(G_{12} \cup G_{13} \cup G_{23})$, the
triangles spanned by these graphs. Define a random $3$-graph $H$
by taking each triangle of $\triangle$ to be an edge with probability
$d_{123}$. With high probability $H$ has density
$d(H) = |E(H)|/n^3 = d_{12}d_{13}d_{23}d_{123} + o(1)$.
Now consider the number of copies of some fixed $3$-graph $F$
in $H$, or for a better parallel with our analysis of quasirandom graphs
consider the `tripartite $F$-homomorphism density'
$d_F(H) = \mb{P}_{\phi \in \Phi}[\phi(F) \sub H]$,
where $\Phi$ is the set of tripartite maps from $F$ to $H$.
We can no longer estimate this by $d(H)^{|E(F)|}$, as might at first be
expected from the analysis for graphs. For a simple example, suppose that
$F$ consists of two edges sharing a pair of vertices; say $F$
has vertex set $A = A_1 \cup A_2 \cup A_3$ with
$A_1 = \{a_1\}$, $A_2 = \{a_2\}$, $A_3 = \{a_3^1,a_3^2\}$
and edges $a_1a_2a_3^1$, $a_1a_2a_3^2$. Then with high probability
$d_F(H) = d_{12}d_{13}^2d_{23}^2d_{123}^2 + o(1)$.

\subsection{Quasirandom complexes}

The second example above shows that a counting lemma for $3$-graphs must
take account of densities of pairs, as well as densities of triples.
Thus we are led to define quasirandomness not for a $3$-graph in isolation,
but for a simplicial complex consisting of a $3$-graph together with all
subsets of its edges. We say that $H$ is a tripartite $3$-complex
\footnote{We prefer this term to `chain', which is used in \cite{G1}}
on $X = X_1 \cup X_2 \cup X_3$ if we have $H = \cup_{I \sub \{1,2,3\}} H_I$,
where $H_\emptyset = \{\emptyset\}$, $H_{\{i\}}$ is a subset of $X_i$ for $1 \le i \le 3$,
$H_{\{i,j\}}$ is a bipartite graph with parts $H_{\{i\}}$, $H_{\{j\}}$ for $1 \le i<j \le 3$,
and $H_{\{1,2,3\}}$ is a $3$-graph contained in the set of triangles spanned
by $H_{\{1,2\}}$, $H_{\{1,3\}}$ and $H_{\{2,3\}}$. Of course, the interesting part
of this structure is the $3$-graph together with its underlying graphs:
we usually take $H_{\{i\}} = X_i$ for $1 \le i \le 3$, and
we only include $H_\emptyset$ so that we are formally correct when referring
to $H$ as a (simplicial) complex. Sometimes we allow $H_{\{i\}}$ to be a strict
subset of $X_i$, but then we can reduce to the usual case by redefining the
ground set as $X' = X'_1 \cup X'_2 \cup X'_3$, where $X'_i = H_{\{i\}}$, $1 \le i \le 3$.

As usual we identify each $H_I$ with its characteristic function, for example
$H_{123}(x_1,x_2,x_3)$ is $1$ if $x_1x_2x_3$ is an edge of $H_{123}$, otherwise $0$
(henceforth we write $123$ instead of $\{1,2,3\}$, etc. for more compact notation).
We let $H_{123}^*$ denote the set of triangles spanned by $H_{\{1,2\}}$, $H_{\{1,3\}}$
and $H_{\{2,3\}}$, and also the characteristic function of this set. Then we define
the relative density by
$$d_{123}(H) = |H_{123}|/|H_{123}^*|.$$
In words, it is the proportion of graph triangles that are triples of the complex.
Note that we can describe the densities of the bipartite graph $H_{ij}$ with a similar
notation: we let $H_{ij}^*$ denote the pairs spanned by $H_i$ and $H_j$, i.e.
the complete bipartite graph with parts $X_i$ and $X_j$, and then $H_{ij}$ has
density $d_{ij}(H)=|H_{ij}|/|H_{ij}^*|$.

We have seen that vertex-regularity is not the correct notion for defining
quasirandomness in $3$-graphs; it turns out that the other two properties described
in Theorem \ref{c4-tfae} {\em can} be generalised to a useful concept for $3$-graphs.
These properties were defined using $4$-cycles for graphs; for $3$-graphs they
are defined in terms of the {\em octahedron} $O_3$. This is a tripartite $3$-graph
with vertex set $A = A_1 \cup A_2 \cup A_3$, where $A_i = \{a_i^0,a_i^1\}$, $1 \le i \le 3$
and $8$ edges $a_1^{e_1}a_2^{e_2}a_3^{e_3}$, where $e=(e_1,e_2,e_3)$ ranges over
vectors in $\{0,1\}^3$. One could define quasirandomness of a complex in terms of the
density of octahedra (by analogy with property 2 in Theorem \ref{c4-tfae}), but we will
follow Gowers and consider `octahedra with respect to the balanced function'
(by analogy with property 3 in Theorem \ref{c4-tfae}).
The {\em balanced function} is defined as
$$\ov{H}_{123}(x_1,x_2,x_3) = H_{123}(x_1,x_2,x_3) - d_{123}(H) H_{123}^*(x_1,x_2,x_3),$$
i.e. for $x_1 \in X_1$, $x_2 \in X_2$, $x_3 \in X_3$,
$\ov{H}_{123}(x_1,x_2,x_3)$ equals $1-d_{123}(H)$ if $x_1x_2x_3$ is
an edge of $H_{123}$, equals $-d_{123}(H)$ if $x_1x_2x_3$ is a triangle in
$H_{12} \cup H_{13} \cup H_{23}$ but not an edge of $H_{123}$, or equals $0$ otherwise.
As before we note that
$\mb{E}_{x_1 \in X_1, x_2 \in X_2, x_3 \in X_3} \ov{H}_{123}(x_1,x_2,x_3) = 0$.
Let $\Phi$ be the set of tripartite maps from $A = A_1 \cup A_2 \cup A_3$
to $X = X_1 \cup X_2 \cup X_3$.
For any function $f:X_1 \times X_2 \times X_3 \to [-1,1]$ we define
$$O_3(f) = \mb{E}_{\phi \in \Phi} \prod_{e \in E(O_3)} f(\phi(e)).$$
We say that $f$ is $\eta$-quasirandom (with respect to $H$)
if $O_3(f) < \eta (d_{12}d_{13}d_{23})^4$
(note that we are assuming that $d_1=d_2=d_3=1$).
One can think of this as saying that $O_3(f)$ is small compared to the density
of the `graph octahedron', i.e the complete tripartite graph $K(2,2,2)$, as
this has density about $(d_{12}d_{13}d_{23})^4$ when the graphs $H_{ij}$ are
quasirandom (as they will be). Then we say that $H_{123}$ is $\eta$-quasirandom
(with respect to $H$) if its balanced function $\ov{H}_{123}$ is $\eta$-quasirandom.

\subsection{Quasirandom decomposition and the counting lemma}

Next we will see why the notion of quasirandomness given in the previous
subsection is useful: there is a decomposition theorem and a counting lemma.

To describe the decomposition theorem we need to first think about the kind of
structure that will arise in a decomposition into simplicial complexes. It will
not simply be a disjoint union of simplicial complexes: just as in a graph
decomposition the various bipartite graphs may share vertices, in a $3$-graph
decomposition the various tripartite $3$-complexes may share vertices and pairs.
Suppose $H$ is an $r$-partite $3$-graph on $X = X_1 \cup \cdots \cup X_r$.
A decomposition of $H$ is described by partitions
$X_i = X_i^1 \cup \cdots \cup X_i^{n_i}$ of each part $X_i$ into subsets,
and partitions $K(X_i,X_j) = G_{ij}^1 \cup \cdots \cup G_{ij}^{n_{ij}}$
of each complete bipartite graph with parts $X_i$, $X_j$ into bipartite subgraphs.
Let $P$ denote this `partition system'. For each edge $e$ of $H$ there is
an {\em induced complex} $H(e,P)$ defined as follows. The vertex set of
$H(e,P)$ is $X_i^{a_i} \cup X_j^{a_j} \cup X_k^{a_k}$, where
$e = x_ix_jx_k$ with $x_i \in X_i^{a_i}$, $x_j \in X_j^{a_j}$, $x_k \in X_k^{a_k}$.
The pairs in $H(e,P)$ are given by the
restriction of $G_{ij}^{a_{ij}} \cup G_{ik}^{a_{ik}} \cup G_{jk}^{a_{jk}}$
to the vertex set of $H(e,P)$,
where $x_ix_j \in G_{ij}^{a_{ij}}$, $x_i x_k \in G_{ik}^{a_{ik}}$, $x_j x_k \in G_{jk}^{a_{jk}}$.
The triples in $H(e,P)$ are those edges of $H$ that also form
triangles in the pairs of $H(e,P)$.

The following decomposition is a variant of Theorem 8.10 in \cite{G1}. (We defer a
justification of the differences here and elsewhere until the appendix.)

\begin{theo} \label{decomp3}
Suppose $0 < 1/n \ll d_1 \ll \eta_2 \ll d_2 \ll \eta_3 \ll 1/r, \eps < 1$
and $H$ is an $r$-partite $3$-graph on $X = X_1 \cup \cdots \cup X_r$,
where each $|X_i| \ge n$.
Then there is a partition system $P$ consisting of
partitions $X_i = X_i^0 \cup X_i^1 \cup \cdots \cup X_i^{n_i}$
with each $n_i \le 1/d_1$
and partitions $K(X_i,X_j) = G_{ij}^1 \cup \cdots \cup G_{ij}^{n_{ij}}$
with each $n_{ij} \le d_2^{-1/2}$,
such that,
\begin{itemize}
\item each exceptional class $X_i^0$ has size at most $\eps|X_i|$, for $1 \le i \le r$,
\item the classes $X_i^{a_i}$, $1 \le i \le r$, $1 \le a_i \le n_i$ all have the same size,
and
\item if $e$ is a randomly chosen edge of $H$, then
with probability at least $1-\eps$, in the induced complex $H(e,P)$,
the graphs are $\eta_2$-quasirandom,
and the $3$-graph is $\eta_3$-quasirandom with respect to $H(e,P)$.
\end{itemize}
\end{theo}

To apply this decomposition we need a counting lemma,
estimating the number of homomorphisms from
a fixed $r$-partite $3$-complex $F$ on $Y = Y_1 \cup \cdots \cup Y_r$
to a quasirandom $r$-partite $3$-complex $H$ on $X = X_1 \cup \cdots \cup X_r$.
Let $\Phi(Y,X)$ denote the set of $r$-partite maps $\phi:Y \to X$.
We define the partite homomorphism density of $F$ in $H$ by
\footnote{Note that we use the same notation for the normal homomorphism density:
it will be clear from the context which is intended, and in any case they are roughly
equivalent for our purposes in that they only differ by a constant factor.}
$$d_F(H) = \mb{P}_{\phi \in \Phi(Y,X)} \prod_{e \in F} H(\phi(e)).$$
The following theorem is a variant of Corollary 5.2 in \cite{G2}.

\begin{theo} \label{count3}
Suppose $0 \le \eta_2 \ll d_2 \ll \eta_3 \ll d_3, \eps, 1/|F| < 1$, that
$H$ and $F$ are $r$-partite $3$-complexes on
$X = X_1 \cup \cdots \cup X_r$ and $Y = Y_1 \cup \cdots \cup Y_r$ respectively,
every graph $H_{ij}$ is $\eta_2$-quasirandom with density at least $d_2$, and
that every $3$-graph $H_{ijk}$ is $\eta_3$-quasirandom with relative density
at least $d_3$ with respect to $H$.
Then
\footnote{Recall that $d_e(H)$ means the density $d_I(H)$ for which $e \in H_I$.}
$$d_F(H) = \mb{P}_{\phi \in \Phi(Y,X)} \prod_{e \in F} H(\phi(e))
= (1 \pm \eps) \prod_{e \in F} d_e(H).$$
\end{theo}

\nib{Remarks.} (1) The special case when $F$ is a single edge shows that
the absolute density of a triple in $H$ is well-approximated by the product
of its relative densities, viz.
\begin{equation} \label{density3}
d(H_{ijk})=|H_{ijk}|/|X_i||X_j||X_k|=(1\pm\eps)d_{ij}d_{ik}d_{jk}d_{ijk}.
\end{equation}

(2) It is important to note the hierarchy of the parameters, as this is where some
of the technical difficulties in hypergraph regularity lie. In the decomposition
the number of graphs $G_{ij}$ in the partition may be much larger than $\eta_3^{-1}$,
so the counting lemma has to cope with graphs $H_{ij}$ with density much smaller
than $\eta_3$, the quasirandomness parameter for triples; fortunately it can
be arranged that $\eta_2$, the quasirandomness parameter for pairs, is smaller still,
and this is sufficient for the counting lemma.

\subsection{Counting homomorphisms to functions and uniform edge-distribution}

We will also need the following more general theorem,
a variant of Theorem 5.1 in \cite{G2}, which allows estimation
of the number of $F$-homomorphisms with respect to functions supported on $H$;
the generalisation is similar to that which we saw earlier from counting
$4$-cycles in quasirandom graphs to estimating $C_4(f)$ for a function $f$.
Intuitively it expresses the fact that, in a quasirandom $3$-complex, sets of size
$k \le 3$ are almost uncorrelated with sets of size less than $k$.
First we need some notation. Suppose $X = X_1 \cup \cdots \cup X_r$ is
a set partitioned into $r$ parts. For $A \sub [r]$ let $K_A(X)$ be the
set of all $A$-tuples with one point in each $X_i$, $i \in A$. We say
that $e \in K_A(X)$ has {\em index} $A$ and we sometimes abuse notation
and use $e$ instead of $A$; thus, if $H$ is an $r$-partite $3$-complex
on $X$ and $e \in H_A$ for some $A \sub [r]$, $|A| \le 3$ then $H_e=H_A$.

\begin{theo} \label{gencount3}
Suppose $0 \le \eta_2 \ll d_2 \ll \eta_3 \ll d_3, \eps, 1/|F| < 1$, that
$H$ and $F$ are $r$-partite $3$-complexes on
$X = X_1 \cup \cdots \cup X_r$ and $Y = Y_1 \cup \cdots \cup Y_r$ respectively,
every graph $H_{ij}$ is $\eta_2$-quasirandom with density at least $d_2$, and
that every $3$-graph $H_{ijk}$ is $\eta_3$-quasirandom with relative density
at least $d_3$ with respect to $H$.
Suppose $F_0$ is a subcomplex of $F$, for each $e \in F$ we have a function
$f_e$ on $K_e(X)$ with $f_e=f_eH_e$,
\footnote{Here we are using the notation described before the theorem:
if $e$ has index some $A \sub [r]$, $|A| \le 3$, then $K_e(X)=K_A(X)=\prod_{i\in A}X_i$.
The equation $f_e=f_eH_e$ is a concise way of saying that $f_e$ is supported
on $H_e=H_A$.}
and $f_e=H_e$ for all $e \in F \sm F_0$.
Then
$$\mb{E}_{\phi \in \Phi(Y,X)} \left[ \prod_{e \in F} f_e(\phi(e)) \right]
= \prod_{e \in F \sm F_0} d_e(H)
\cdot \mb{E}_{\phi \in \Phi(Y,X)} \left[ \prod_{e \in F_0} f_e(\phi(e)) \right]
\pm \eps \prod_{e \in F} d_e(H).$$
\end{theo}

Note that Theorem \ref{count3} is the case of Theorem \ref{gencount3} when
$F_0=\emptyset$.
The following corollary shows that quasirandomness implies vertex-uniformity,
the analogue of graph regularity that we mentioned earlier as being too weak
for a useful hypergraph regularity theory.

\begin{coro} \label{vertex-uniform3}
Suppose $0 \le \eta_2 \ll d_2 \ll \eta_3 \ll d_3, \eps < 1$, that
$H$ is a tripartite $3$-complex on $X = X_1 \cup X_2 \cup X_3$
with $H_i=X_i$, $1 \le i \le 3$,
that $H_{12}$, $H_{13}$ and $H_{23}$ are $\eta_2$-quasirandom with density at least $d_2$,
and that $H_{123}$ is $\eta_3$-quasirandom with relative density
at least $d_3$ with respect to $H$.
Suppose $W_i \sub X_i$, $1 \le i \le 3$ and let $H_{123}[W]$ be the restriction
of $H_{123}$ to $W = W_1 \cup W_2 \cup W_3$.
Then
$$\frac{|H_{123}[W]|}{|X_1||X_2||X_3|} = \frac{|H_{123}|}{|X_1||X_2||X_3|} \cdot
\left( \frac{|W_1||W_2||W_3|}{|X_1||X_2||X_3|} \pm \eps \right).$$
\end{coro}

\nib{Proof.}
Set $Y_i = \{i\}$ for $1 \le i \le 3$, $F$ equal to all subsets
of $Y = Y_1 \cup Y_2 \cup Y_3$ and $F_0 \sub F$ equal to the
subsets of size at most $1$. Apply Theorem \ref{gencount3} with
$f_i(x_i)=1_{x_i \in W_i}$ and $f_e=H_e$ for $e \in F \sm F_0$.
Since $H_i=X_i$, $1 \le i \le 3$
we have $\prod_{e \in F \sm F_0} d_e(H) = \prod_{e \in F} d_e(H)$.
Also
$\mb{E}_{\phi \in \Phi(Y,X)} \prod_{e \in F_0} f_e(\phi(e))
= \prod_{i=1}^3 |W_i|/|X_i|$.
Now
$\mb{E}_{\phi \in \Phi(Y,X)} \left[ \prod_{e \in F} f_e(\phi(e)) \right]
= |H_{123}[W]|/|X_1||X_2||X_3|$
so Theorem  \ref{gencount3} gives the result. \qed

\subsection{Uniform distribution of octahedra in quasirandom $3$-complexes}

Now we will establish a further special case of our quasirandom counting lemma,
the uniform distribution of octahedra in quasirandom $3$-complexes,
which is the heart of our solution of Mubayi's conjecture.
Recall that the octahedron $O_3$ is a tripartite $3$-graph
with vertex set $A = A_1 \cup A_2 \cup A_3$,
where $A_i = \{a_i^0,a_i^1\}$, $1 \le i \le 3$
and $8$ edges $a_1^{e_1}a_2^{e_2}a_3^{e_3}$, where $e=(e_1,e_2,e_3)$ ranges over
vectors in $\{0,1\}^3$. We write $O_3^\le$ for the tripartite $3$-complex
generated by $O_3$ (i.e. it contains all subsets of the edges of $O_3$)
and $O_3^< = O_3^{\le} \sm O_3$ for the tripartite $2$-complex of strict
subsets of edges of $O_3$.

Suppose $H$ is a tripartite $3$-complex on $X = X_1 \cup X_2 \cup X_3$,
where $H_i=X_i$, $1 \le i \le 3$. We define an `auxiliary'
tripartite $3$-complex $H'$ on $X' = X'_1 \cup X'_2 \cup X'_3$,
where $X'_i$ consists of all maps $\phi_i:A_i \to X_i$, $1 \le i \le 3$
as follows. The triples $H'_{123}$ consist of all
$\phi_1\phi_2\phi_3$ for which $\phi=(\phi_1,\phi_2,\phi_3)$
is a homomorphism from $O_3$ to $H$;
the pairs $H'_{ij}$ consist of all
$\phi_i\phi_j$ such that $(\phi_i,\phi_j)$ is a homomorphism
from $(O_3)^\le_{ij}$ to $H_{ij}$, for $1 \le i<j \le 3$.

\begin{theo} \label{qr-oct}
Suppose $0 \le 1/n \ll \eta_2 \ll \eps_2 \ll \eta'_2 \ll d_2
\ll \eta_3 \ll \eps_3 \ll \eta'_3 \ll \eps', d_3 < 1$, that
$H$ is a tripartite $3$-complex on $X = X_1 \cup X_2 \cup X_3$,
where $H_i=X_i$ with $|X_i| \ge n$, $1 \le i \le 3$,
that $H_{12}$, $H_{13}$ and $H_{23}$ are $\eta_2$-quasirandom with density at least $d_2$,
and that $H_{123}$ is $\eta_3$-quasirandom with relative density
at least $d_3$ with respect to $H$. Then, in the auxiliary complex $H'$,
each $H'_{ij}$ is $\eta'_2$-quasirandom with density $(1 \pm \eps_2)d_{ij}(H)^4$,
and $H'_{123}$ is $\eta'_3$-quasirandom with relative density $(1 \pm \eps_3)d_{123}(H)^8$
with respect to $H'$.

Furthermore, suppose that $G_i$ is a graph on $X_i$ with $|E(G_i)|> \eps'|X_i|^2$, $1 \le i \le 3$
and we randomly and independently select edges $x_ix'_i$ of $G_i$, $1 \le i \le 3$.
Then $\{x_1x'_1,x_2x'_2,x_3x'_3\}$ spans an octahedron in $H$ with probability
$(1 \pm \eps') d_{12}(H)^4d_{13}(H)^4d_{23}(H)^4d_{123}(H)^8$.
\end{theo}

\nib{Proof.} It may help the reader to note that the prime symbol $'$ is a visual
cue, distinguishing quantities referring to $H'$ from those referring to $H$.
We write $d_I=d_I(H)$ and $d'_I=d_I(H')$ for $I \sub \{1,2,3\}$, $|I|=2,3$.
\footnote{To avoid confusion, note that $d_2$, $d_3$ are lower bounds for
$d_{ij}$, $d_{123}$ respectively, not relative densities of $H_2$, $H_3$.}
The quasirandomness of each $H'_{ij}$ follows from the discussion of quasirandom counting
of $4$-cycles in subsection 2.3, so it remains to establish the claim for $H'_{123}$.
The final statement will then follow from Corollary \ref{vertex-uniform3}
and equation (\ref{density3}).

First we estimate the relative density of $H'$ using equation (\ref{density3}):
we have (using $\eps_3/3$ instead of $\eps$)
$$d(H'_{123}) =|H'_{123}|/|X_1|^2|X_2|^2|X_3|^2
=\mb{P}_{\phi \in \Phi(A,X)}[\phi(O_3^\le) \sub H]
= (1 \pm \eps_3/3) d_{12}^4d_{13}^4d_{23}^4d_{123}^8,$$
and
$$d(H^{\prime *}_{123}) =|H^{\prime *}_{123}|/|X_1|^2|X_2|^2|X_3|^2
=\mb{P}_{\phi \in \Phi(A,X)}[\phi(O_3^<) \sub H]
= (1 \pm \eps_3/3) d_{12}^4d_{13}^4d_{23}^4,$$
so $d_{123}(H')=|H'_{123}|/|H^{\prime *}_{123}|=(1 \pm \eps_3)d_{123}^8$.

\begin{figure}
\begin{center}
\includegraphics[height=7cm]{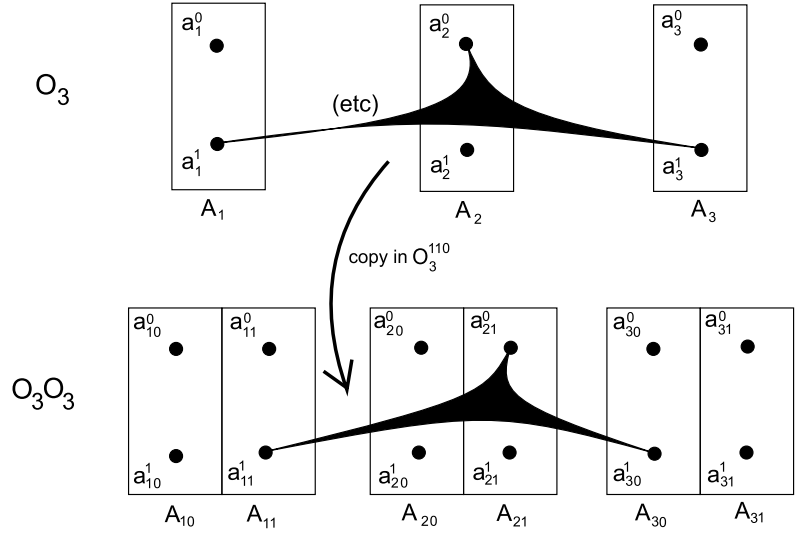}
\end{center}
\caption{$O_3$ and $O_3O_3$: the left index identifies part $1$, $2$ or $3$; the right index
identifies the copy $O_3^e$ of $O_3$; the top index determines which edge we consider within
a copy of $O_3$.}
\label{oct-pic}
\end{figure}

Next we need to estimate $O_3(\ov{H'}_{123})$. Write
$f'_0 = H'_{123}$ (the characteristic function) and
$f'_1 = d'_{123}H^{\prime *}_{123}$,
so that $\ov{H'}_{123} = f'_0-f'_1$.
Then
$$O_3(\ov{H'}_{123}) = \mb{E}_{\phi' \in \Phi(A,X')} \prod_{e \in O_3} \ov{H'}_{123}(\phi'(e))
= \sum_j (-1)^{\sum j} \mb{E}_{\phi' \in \Phi(A,X')} \prod_{e \in O_3} f'_{j_e}(\phi'(e)),$$
where $j$ ranges over all $\{0,1\}$-sequences $(j_e)_{e \in O_3}$ and
$\sum_j = \sum_{e \in O_3} j_e$.
We interpret each summand as counting homomorphisms to functions from $O_3O_3$,
by which we mean (with unusual notation!) the complete tripartite $3$-graph
with $4$ vertices in each class (cf. our earlier analysis of $C_4$ via $K_{4,4}$.)
To accomplish this we consider a complete tripartite $3$-graph, which we call $O_3O_3$,
with $3$ parts $A_{i,0} \cup A_{i,1}$, $1 \le i \le 3$,
where $A_{i,j} = \{a_{i,j}^0,a_{i,j}^1\}$, $1 \le i \le 3$, $j=0,1$.
We think of the edges of $O_3O_3$ as a union of $8$ copies of $O_3$, indexed by $O_3$,
thus $E(O_3O_3) = \cup_{e \in O_3} O_3^e$, where, writing
$e=\{a_1^{e_1},a_2^{e_2},a_3^{e_3}\}$, $O_3^e$ is the copy of $O_3$ with
parts $\{a_{i,e_i}^0,a_{i,e_i}^1\}$, $1 \le i \le 3$ (see figure \ref{oct-pic}).
\footnote{Some readers may like to think of this construction as a `tripartite tensor product'.}

We identify a tripartite
map $\phi'$ from $A = A_1 \cup A_2 \cup A_3$ to $X' = X'_1 \cup X'_2 \cup X'_3$
with a tripartite map $\phi$ from $AA = (A_{1,0} \cup A_{1,1}) \cup (A_{2,0} \cup A_{2,1})
 \cup (A_{3,0} \cup A_{3,1})$ to $X$ by defining
$\phi(a_{i,j}^k) = \phi'(a_i^j)(a_i^k)$. Thus for
$e = \{a_1^{e_1},a_2^{e_2},a_3^{e_3}\} \in O_3$,
$\phi^e=(\phi'(a_1^{e_1}),\phi'(a_2^{e_2}),\phi'(a_1^{e_2}))$
acts on $O_3^e$.
We can rewrite the summands above using
$$\prod_{e \in O_3} f'_{j_e}(\phi'(e))
= \prod_{e \in O_3} \prod_{\hat{e} \in O_3^e} f_{j_e}(\phi^e(\hat{e}))
= \prod_{E \in O_3O_3} f_{j_{e(E)}}(\phi(E)),$$
where we write $e(E)$ for that $e \in O_3$ such that $E \in O_3^e$,
$f_0 = H_{123}$ (the characteristic function) and
$f_1 = (d'_{123})^{1/8}H^*_{123} = (1 \pm \eps_3/4)d_{123}H^*_{123}$.

Now we can estimate the contribution from the summand corresponding to
$j=(j_e)_{e \in O_3}$ using the counting lemma for
homomorphisms to $H$ from a tripartite $3$-complex $K_j$ on $AA$, in which the
pairs form the complete tripartite graph on $AA$ and the triples are
$\cup_{e:j_e=0} O_3^e$. We have
\begin{align*}
\mb{E}_{\phi' \in \Phi(A,X')} \prod_{e \in O_3} f'_{j_e}(\phi'(e))
& = \mb{E}_{\phi \in \Phi(AA,X)} \prod_{E \in O_3O_3} f_{j_{e(E)}}(\phi(E)) \\
& = (d'_{123})^{\sum j} \mb{P}_{\phi \in \Phi(AA,X)}[\phi(K_j) \sub H] \\
& = ((1 \pm \eps_3)d_{123}^8)^{\sum j}
\cdot (1 \pm \eps_3) d_{12}^{16}d_{13}^{16}d_{23}^{16}d_{123}^{64-8\sum j} \\
& = (1 \pm O(\eps_3)) (d'_{12} d'_{13} d'_{23})^4 (d'_{123})^8,
\end{align*}
where the implicit constant in the $O(\cdot)$ notation is absolute, say $100$.
Therefore
\begin{align*}
O_3(\ov{H'}_{123}) & = \sum_j (-1)^{\sum j}
(1 \pm O(\eps_3)) (d'_{12} d'_{13} d'_{23})^4 (d'_{123})^8 \\
& = \sum_j O(\eps_3) (d'_{12} d'_{13} d'_{23})^4 (d'_{123})^8
< \eta'_3 (d'_{12} d'_{13} d'_{23})^4,
\end{align*}
since $\eps_3 \ll \eta'_3$.
This proves that $H'_{123}$ is $\eta'_3$-quasirandom with relative density
$(1 \pm \eps)d_{123}(H)^8$ with respect to $H$.
The third statement of the theorem follows from Corollary \ref{vertex-uniform3},
so we are done. \qed

\section{The Fano plane}

In this section we prove Theorem \ref{fano}, thus establishing Mubayi's conjecture.
Our key tool will be Theorem \ref{qr-oct}, our result above on uniform distribution
of octahedra in quasirandom $3$-complexes. We also make use of the idea of stability,
or approximate structure, which can be traced back to work of Erd\H{o}s and
Simonovits in the 60's in extremal graph theory. Informally stated,
a stability result tells us about the structure of configurations
that are close to optimal in an extremal problem: for example, a
triangle-free graph with $n^2/4 - o(n^2)$ edges differs from a
complete bipartite graph by $o(n^2)$ edges. Such a result is
interesting in its own right, but somewhat surprisingly it is often
a useful stepping stone in proving an exact result. Indeed, it was
developed by Erd\H os and Simonovits to determine the exact Tur\'an
number for $k$-critical graphs.

\begin{figure}
\begin{center}
\includegraphics[height=6cm]{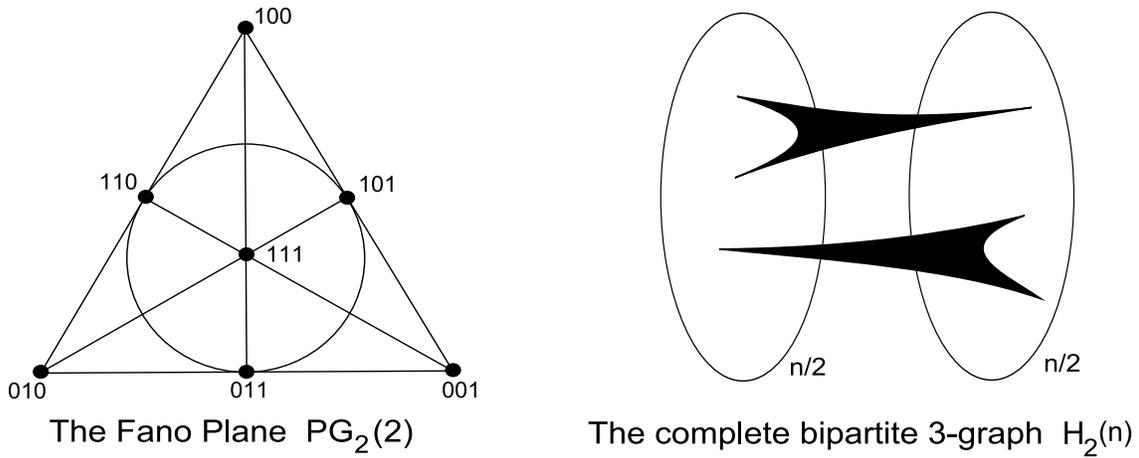}
\end{center}
\caption{$PG_2(2)$ and $H_2(n)$.}
\label{fano-pic}
\end{figure}

To further explain our method we first need some definitions. The {\em Fano plane}
$PG_2(2)$ is the projective plane over $\mb{F}_2 = \{0,1\}$, the field with two elements:
it is a $3$-graph in which the vertex set consists of the $7$ non-zero
vectors of length $3$ over $\mb{F}_2$, and the $7$ edges are those triples of vectors
$abc$ with $a+b=c$ (see figure \ref{fano-pic}).
A natural construction of a $3$-graph not containing the Fano plane
is a {\em balanced complete bipartite $3$-graph} $H_2(n)$, which we recall is
obtained by partitioning a set into two parts
$X = X_1 \cup X_2$ with sizes as equal as possible
and taking as edges all triples that are not contained in
either part (i.e. that intersect both parts): it is easy to verify that the Fano plane
is not bipartite, so is not contained in $H_2(n)$.

In 1976 S\'os \cite{So} conjectured that $H_2(n)$
gives the exact value of $\mbox{ex}(n,PG_2(2))$. This was established for large $n$,
together with the characterisation of $H_2(n)$ as the unique maximising configuration,
independently and simultaneously by Keevash and Sudakov \cite{KS1}
and F\"uredi and Simonovits \cite{FS}.
In \cite{KS1} the following stability method was used:
the first step was to show that a Fano-free $3$-graph with $(3/4-o(1))\binom{n}{3}$
edges differs from a complete bipartite $3$-graph by $o(n^3)$ edges;
then, the second step was to examine the possible imperfections in structure
in a $3$-graph that is close to being bipartite, showing that they exclude more
edges in a Fano-free $3$-graph than a complete bipartite $3$-graph with no imperfections.

In considering the codegree problem, Mubayi \cite{M1} noted that the same construction
$H_2(n)$ gives a lower bound $\mbox{ex}_2(n,PG_2(2)) \ge \lfloor n/2 \rfloor$
and proved that $\pi_2(PG_2(2)) = 1/2$, so this bound is asymptotically best possible.
In our solution of the exact problem, the first step will be to use Theorem \ref{qr-oct}
to obtain some partial structure of Fano $3$-graphs with nearly maximal minimum codegree,
namely that a Fano $3$-graph $H$ on $n$ vertices with minimum $2$-degree $(1/2-o(1))n$
has a sparse set of size $(1/2-o(1))n$. The second step (fairly straightforward)
is to complete the approximate structure, showing that $H$ differs from a complete bipartite
$3$-graph by $o(n^3)$ edges. Then the final step is to examine the possible imperfections
in structure to show that the largest minimum codegree is attained with no imperfections:
this is similar in spirit to the analysis of the normal Tur\'an problem, although there
are some additional technical difficulties for the codegree problem.
We will divide the three steps of the proof into subsections of this section.

\subsection{The sparse set}

The first step of the proof is the following application of Theorem \ref{qr-oct}.

\begin{theo} \label{sparse3}
Suppose $0 < 1/n \ll \theta \ll \eps \ll 1$, $H$ is a $3$-graph on
a set $X$ of $n$ vertices that does not contain a Fano plane, and $\delta_2(H)>(1/2-\theta)n$.
Then there is a set $Z \sub X$ with $|Z|>(1/2-\eps)n$ such that $d(H[Z]) < \eps$,
i.e. $Z$ spans at most $\eps|Z|^3/6$ edges of $H$.
\end{theo}

The idea of the proof is as follows. First we apply Theorem \ref{decomp3} to
decompose $H$ into induced $3$-complexes, most of which are quasirandom.
Then a simple double-counting
gives us a vertex $x$ such that there are about half of the vertex classes
in the partition on which $x$ has a dense neighbourhood graph. For each triple
of such classes, each quasirandom induced $3$-complex on these classes must
be very sparse; otherwise Theorem \ref{qr-oct} gives us an octahedron in
which the pairs from each class are in the neighbourhood graph of $x$,
and this structure contains a Fano plane. This implies that we have about half
of the vertex classes in the partition inducing a sparse subhypergraph of $H$.

\medskip

\nib{Proof of Theorem \ref{sparse3}.}
We operate with a constant hierarchy
$0 < 1/n \ll d_1 \ll \eta_2 \ll d_2 \ll \eta_3 \ll d_3, \gamma \ll \theta \ll \eps \ll 1$.
Set $r=\lceil \theta^{-1} \rceil$ and let $X = X_1 \cup \cdots \cup X_r$ be an arbitrary partition
with $||X_i| - n/r| < 1$ for $1 \le i \le r$.
Let $H_1$ be the edges of $H$ that respect this partition, i.e.
have at most one point in each $X_i$.
Applying Theorem \ref{decomp3} to the generated $r$-partite $3$-complex $H_1^\le$,
using $\gamma$ instead of $\eps$,
we obtain a partition system $P$ consisting of
partitions $X_i = X_i^0 \cup X_i^1 \cup \cdots \cup X_i^{n_i}$
with each $n_i \le 1/d_1$
and partitions $K(X_i,X_j) = G_{ij}^1 \cup \cdots \cup G_{ij}^{n_{ij}}$
with each $n_{ij} \le d_2^{-1/2}$,
such that,
\begin{itemize}
\item each exceptional class $X_i^0$ has size at most $\gamma|X_i|$, for $1 \le i \le r$,
\item the classes $X_i^{a_i}$, $1 \le i \le r$, $1 \le a_i \le n_i$ have the same size, say $m$, and
\item if $e$ is a randomly chosen edge of $H_1$, then
with probability at least $1-\gamma$, in the induced complex $H_1^\le(e,P)$,
the graphs are $\eta_2$-quasirandom,
and the $3$-graph is $\eta_3$-quasirandom with respect to $H_1^\le(e,P)$.
\end{itemize}

Next we find a vertex $x$ as described in the sketch before the proof.
Consider the set of all pairs contained in the vertex classes of $P$:
denote it $Y = \cup_{(i,t_i) \in I} Y_{i,t_i}$,
where $I=\{(i,t_i): 1 \le i \le r, 1 \le t_i \le n_i\}$
and $Y_{i,t_i} = \binom{X_i^{t_i}}{2}$.
(Note that we are not considering any exceptional classes here.)
Double-counting gives
$$\sum_{x \in X} |N_H(x) \cap Y| = \sum_{y \in Y} |N_H(Y) \cap X| > |Y|(1/2-\theta)n,$$
so we can choose $x$ with $|N_H(x) \cap Y|>(1/2-\theta)|Y|$.
Let $C = \{(i,t_i): |N_H(x) \cap Y_{i,t_i}| > \theta |Y_{i,t_i}|$.
Then, since $|X_i^{t_i}|=m$ for every $(i,t_i) \in I$,
$$(1/2-\theta)|I|\binom{m}{2} = (1/2-\theta)|Y| < |N_H(x) \cap Y| < |C| \binom{m}{2}
+ |I \sm C| \theta \binom{m}{2},$$
so $|C| > (1/2-2\theta)|I|$. Let $Z = \cup_{(i,t_i) \in C} X_i^{t_i}$.
Then $|Z| = |C|m > (1/2-3\theta)n$.

Consider any induced tripartite $3$-complex $H_1^\le(e,P)$ in which the three
vertex classes are $X_{i_a}^{t_a}$, $1 \le a \le 3$ with $(i_a,t_a) \in C$, $1 \le a \le 3$
and $i_1,i_2,i_3$ are distinct. We claim that it cannot be that
all three graphs in $H_1^\le(e,P)$ are $\eta_2$-quasirandom with density at least $d_2$
and the $3$-graph in $H_1^\le(e,P)$ is $\eta_3$-quasirandom with density at least $d_3$.
For then, since $\gamma \ll \theta$,
Theorem \ref{qr-oct} gives us an octahedron (in fact many!) with parts
$\{x_a,x'_a\}$, $1 \le a \le 3$ in which $x_ax'_a$ is an edge of $N_H(x) \cap Y_{i_a,t_a}$,
$1 \le a \le 3$. However such a structure clearly contains a Fano plane, which contradicts
our assumptions.

It follows that at least one quasirandomness condition or density condition
fails for each induced tripartite $3$-complex $H_1^\le(e,P)$ with vertex classes in $Z$.
Now the choice of partition system $P$ guarantees that at most $\gamma n^3$
edges of $H_1$ belong to induced complexes which fail a quasirandomness condition.
Also, at most $d_3 n^3$ edges of $H_1$ belong to an induced complex
for which the $3$-graph has relative density at most $d_3$,
and at most $3d_2^{1/2} n^3$ edges of $H_1$ belong to an induced complex
in which one of the $2$-graphs has relative density at most $d_2$
(since $n_{ij} \le d_2^{-1/2}$ for every $i,j$).
Finally, we can estimate the number of edges of $H$ that do not respect
the partition $X = X_1 \cup \cdots \cup X_r$ as $|H \sm H_1| < n^3/r$.
We conclude that $d(H[Z]) < \eps$, as required. \qed

\subsection{Stability}

The remainder of the proof makes no further use of quasirandomness.
The next step is to use the sparse set to deduce a stability result, i.e. that
a Fano-free $3$-graph with minimum codegree about $n/2$ is approximately bipartite.

\begin{theo} \label{stab3}
Suppose $0 < 1/n \ll \theta \ll \eps \ll 1$, $H$ is a $3$-graph on
a set $X$ of $n$ vertices that does not contain a Fano plane, and $\delta_2(H)>(1/2-\theta)n$.
Then there is a partition $X = A \cup B$
such that at most $\eps n^3$ edges are contained entirely within either $A$ or $B$.
\end{theo}

\nib{Proof.}
It is convenient to introduce a constant $\phi$ with $\theta \ll \phi \ll \eps$.
By Theorem \ref{sparse3} we can find $A \sub X$ with $|A|>(1/2-\phi)n$
and $d(H[A]) < \phi$, i.e. $|H[A]| < \phi|A|^3/6$.
Let $B = X \sm A$. Then
$$3|H[A]| + \sum_{p \in \binom{A}{2}} |N_H(p) \cap B|
= \sum_{p \in \binom{A}{2}} |N_H(p)|
> \binom{|A|}{2} (1/2-\theta)n,$$
so
$$\sum_{b \in B} \left|N_H(b) \cap \binom{A}{2}\right|
= \sum_{p \in \binom{A}{2}} |N_H(p) \cap B|
> (1/2-\theta)n \binom{|A|}{2} - \phi|A|^3/2
> (1/2-2\phi)n \binom{|A|}{2}.$$
This implies that $|B|>(1/2-2\phi)n$, so both $|A|$ and $|B|$ are $(1/2 \pm 2\phi)n$.
Let
$$B_0 = \left\{b \in B: \left|N_H(b) \cap \binom{A}{2}\right|
< (1 - \phi^{1/2}) \binom{|A|}{2} \right\}.$$
Then
\begin{align*}
(1/2-2\phi)n \binom{|A|}{2}
& < \sum_{b \in B} \left|N_H(b) \cap \binom{A}{2}\right| \\
& < |B_0| (1 - \phi^{1/2}) \binom{|A|}{2}  + |B \sm B_0| \binom{|A|}{2} \\
& = |B| \binom{|A|}{2} - |B_0| \cdot \phi^{1/2} \binom{|A|}{2},
\end{align*}
so $|B_0| < \phi^{-1/2}(|B|-(1/2-2\phi)n) < 3\phi^{1/2}n$.

We claim that there is no edge $e \in H[B \sm B_0]$. Otherwise we would have
$$\left| \bigcap_{b \in e} N_H(b) \cap \binom{A}{2} \right| > (1-3\phi^{1/2})\binom{|A|}{2}.$$
Then, since $\phi \ll 1$,
$\cap_{b \in e} N_H(b) \cap \binom{A}{2}$ contains a copy of $K_4$,
the complete graph on $4$ vertices: we can see this by
Tur\'an's theorem, which gives $\pi(K_4)=2/3$,
or even just a simple averaging argument which gives $\pi(K_4) \le 5/6$.
However, this $K_4$ and $e$ span a copy of the Fano plane, contradiction.
Therefore every edge of $B$ contains a point of $B_0$, so
$|H[B]| < |B_0||B|^2 < 3\phi^{1/2}n^3$. Since $|H[A]| < \phi|A|^3/6$
there are less than $\eps n^3$ edges contained entirely within either $A$ or $B$. \qed

\subsection{The exact codegree result for the Fano plane}

Finally we use the stability result above to prove Theorem \ref{fano}, which states that
if $n$ is sufficiently large and $H$ is a $3$-graph on $n$ vertices
with minimum $2$-degree at least $n/2$ that does not contain a Fano
plane, then $n$ is even and $H$ is a balanced complete bipartite $3$-graph.

\medskip

\nib{Proof of Theorem \ref{fano}.}
Let $\eps$ be sufficiently small and $n > n_0(\eps)$ sufficiently large.
By Theorem \ref{stab3} we have a partition $X = X_0 \cup X_1$
so that $|H[X_0]|+|H[X_1]| < \eps n^3$. Choose the partition
that minimises $|H[X_0]|+|H[X_1]|$. We will show that
this partition satisfies the conclusion of the theorem.
Note first that the same argument used in the proof of Theorem \ref{stab3}
shows that $|X_0|$ and $|X_1|$ are $(1/2 \pm 2\eps)n$.

Next we show that there is no vertex $x \in X_0$ with degree at least
$\eps^{1/4} n^3$ in $H[X_0]$. For suppose there is such a vertex $x$.
By choice of partition we have $|N_H(x) \cap \binom{X_1}{2}|
\ge|N_H(x) \cap \binom{X_0}{2}| \ge \eps^{1/4}n^2$, or we could
reduce $|H[X_0]|+|H[X_1]|$ by moving $x$ to $X_1$.
We can choose matchings
$M = \{x^1_1 x^2_1,\cdots, x^1_m x^2_m\}$
in $N_H(x) \cap \binom{X_0}{2}$
and $M' = \{y^1_1 y^2_1,\cdots,y^1_m y^2_m\}$
in $N_H(x) \cap \binom{X_1}{2}$,
with $m = \eps^{1/4} n/2$;
indeed, it is a well-known observation that any maximal matchings
will be at least this large, as the vertex set of a maximal matching
in a graph covers all of its edges.
Now we use an averaging argument to find a pair of edges of $M$
that are not `traversed' by many edges of $H$, in the following sense.
We have
$$ \sum_{I = \{i_1,i_2\} \in \binom{[m]}{2}}
\sum_{J = \{j_1,j_2\} \in \{1,2\}^2}
 |N_{H[X_0]}(\{x^{j_1}_{i_1},x^{j_2}_{i_2}\})|
< 3|H[X_0]| < 3\eps n^3,$$
so we can choose $I \in \binom{[m]}{2}$ such that
$$\sum_{J = \{j_1,j_2\} \in [2]^2}
|N_{H[X_0]}(\{x^{j_1}_{i_1},x^{j_2}_{i_2}\})|
< \binom{m}{2}^{-1} 3\eps n^3 < 50 \eps^{1/2} n.$$
Also, for every $1 \le k  \le m$, $H$ cannot contain
all $8$ edges $x^{j_1}_{i_1} x^{j_2}_{i_2} y_k^{j_3}$, with
$j_1,j_2,j_3$ in $\{1,2\}$, as then together with $x$
we would have a copy of the Fano plane.
By double-counting we deduce that
there is some pair $p = x^{j_1}_{i_1} x^{j_2}_{i_2}$
such that there are at least $m/4$ vertices $y_k^{j_3}$ that
do not belong to $N_H(p)$.
This gives
$|N_H(p)| < 50\eps^{1/2} n + |X_1| - \eps^{1/4} n/8 < (1/2-\eps)n$,
which contradicts our assumptions. We deduce that there
is no vertex in $X_0$ with degree at least $\eps^{1/4} n^2$ in $H[X_0]$.
Similarly there is no vertex in $X_1$
with degree at least $\eps^{1/4} n^2$ in $H[X_1]$.

Write $|X_0|=n/2+t$ and $|X_1|=n/2-t$, where without loss of
generality $0 \le t \le 2\eps n$. Suppose for a contradiction that either
$t>0$ or $t=0$ and there is an edge in $H[X_0]$ or $H[X_1]$.
Note that for every pair $p \in \binom{X_0}{2}$
we have $|N_H(p) \cap X_0| \ge |N_H(p)| - |X_1| \ge t$.
Thus we can assume there is at least one edge in $X_0$ (since the case
$t=0$ is symmetrical).
Let $Y_i \sub X_i$ be minimum size
transversals of $H[X_i]$ (i.e. $|e \cap Y_i| \ne \emptyset$ for
every $e$ in $H[X_i]$). Then $Y_0 \ne \emptyset$.
Also, by the previous paragraph
$$|Y_0|\eps^{1/4}n^2 > \sum_{x \in Y_0} |N_{H[X_0]}(x)| \ge |H[X_0]|
= \frac{1}{3} \sum_{p \in \binom{X_0}{2} } |N_{H[X_0]}(p)|
\ge \frac{t}{3} \binom{|X_0|}{2},$$
so $|Y_0| > \eps^{-1/4} \frac{t}{3} n^{-2} \binom{n/2+t}{2}
> 50t$, say, since $\eps$ is small.

For each edge $e \in X_0$ consider all possible ways to extend it
to a copy of the Fano plane using some $F \in \binom{X_1}{4}$.
Since $H$ does not contain a Fano plane there is some triple
with $2$ points in $F$ and one point in $e$ that is not an edge.
We count each such triple at most $\binom{|X_1|-2}{2}$ times, so
we get a set of at least $\binom{|X_1|-2}{2}^{-1} \binom{|X_1|}{4}$
distinct triples.
Thus there is some point $x$ in $e$ for which we have
a set $M_x$ of at least
$\frac{1}{3} \binom{|X_1|-2}{2}^{-1} \binom{|X_1|}{4}
> \frac{1}{20} \binom{|X_1|}{2}$ `missing' triples involving
$x$ and a pair in $\binom{X_1}{2}$ not belonging to $N_H(x)$.
Now varying $e$ over all edges in $X_0$, we get sets $M_x$, $x \in T$
of at least $\frac{1}{20} \binom{|X_1|}{2}$ `missing' triples involving $x$,
for some transversal $T$. Since $Y_0$ is a minimum size transversal,
$M = \cup_{x \in T} M_x$ has size at least
$|M| \ge |Y_0| \cdot \frac{1}{20} \binom{|X_1|}{2}$.

Now, for each pair $p \in \binom{X_1}{2}$ we have
$$n/2 \le |N_H(p)|  = |N_H(p) \cap X_0| + |N_H(p) \cap X_1|
\le |X_0|-|N_M(p)| + |N_H(p) \cap X_1|,$$
so
$$3|H[X_1]| = \sum_{p \in \binom{X_1}{p}} |N_H(p) \cap X_1|
 \ge \sum_{p \in \binom{X_1}{2}} (|N_M(Q)|-t)
 = |M| - t\binom{|X_1|}{2}.$$
Since $|M| \ge  \frac{1}{20}|Y_0| \binom{|X_1|}{2}$ and
$\frac{1}{20}|Y_0| - t > \frac{1}{40}|Y_0|$ we have
$|H[X_1]| > \frac{1}{120} |Y_0|\binom{|X_1|}{2}$.
Also $|Y_1|\eps^{1/4}n^2 > \sum_{x \in Y_1} |N_{H[Y_1]}(x)| > |H[X_1]|$,
so $|Y_1| >  \eps^{-1/4}  \frac{1}{120} |Y_0| n^{-2} \binom{n/2-t}{2}
> 2|Y_0|$, say, since $\eps$ is small.

Finally we can apply the argument of the previous paragraph
interchanging $X_0$ and $X_1$. We get a set $M'$ of at
least $|Y_1| \cdot \frac{1}{20} \binom{|X_0|}{2}$ distinct triples
that are not edges, each having $2$ points in $X_0$ and
$1$ point in $X_1$. For each pair $p \sub X_0$ we now have
$|N_H(p) \cap X_0| \ge |N_{M'}(p)| + t$, so
$|H[X_0]| > \frac{1}{3} \cdot (|Y_1|/20+t)\binom{|X_0|}{2}$
and
$|Y_0| > \eps^{-1/4}n^{-2} |H[X_0]| > 2|Y_1|$.
This contradiction completes the proof. \qed

\section{Quasirandom hypergraphs}

We now return to the theory of quasirandom hypergraphs: we will
discuss the Gowers approach in full generality, together with
some variants that we need for our arguments, the general
form of our quasirandom counting lemma, and its application to
generalised Tur\'an problems.
The essential ideas are already present in the discussion above
for graphs and $3$-graphs, so we will be fairly brief.
We refer the reader to \cite{G2} for full details of
the Gowers theory, with the proviso that we have also adopted some
notation and terminology from \cite{GT} and made
other changes for consistency of notation.
This section is divided into four subsections: the first extends the
notation and definitions introduced above to general $k$-graphs,
the second contains the general forms of the decomposition theorem
and counting lemma, the third the general form of our quasirandom
counting lemma, and the fourth a generalised form of Theorem \ref{sparse3}.

\subsection{Definitions}

Consider a set $X = X_1 \cup \cdots \cup X_r$ partitioned into
$r$ parts. For brevity call this an $r$-partite set.
$A \sub X$ is an $r$-partite subset if it contains at most
one point from each $X_i$. An $r$-partite $k$-graph on $X$ is
a $k$-graph on $X$ consisting of $r$-partite subsets. An $r$-partite
$k$-complex $H$ is a collection of $r$-partite subsets of $X$ of
size at most $k$ that forms a simplicial complex, i.e. if
$S \in H$ and $T \sub S$ then $T \in H$. If $H$ is a $k$-graph
then we can generate a $k$-complex
$H^{\le} = \{T: \exists S \in H, T \sub S\}$ and a $(k-1)$-complex
$H^{<} = \{T: \exists S \in H, T \subset S\}$ (strict subsets).
If $H$ is an $r$-partite $k$-complex on $X = X_1 \cup \cdots \cup X_r$
and $A \in \binom{[r]}{\le k}$ we define an $|A|$-complex
$H_{A^\le} = \cup_{B \sub A} H_B$ and an $(|A|-1)$-complex
$H_{A^<} = \cup_{B \subset A} H_B$.

The index of an $r$-partite subset $A$ of $X$ is
$i(A) = \{ i \in [r]: A \cap X_i \neq \emptyset \}$.
Let $K_A(X) = \{S \sub X: i(S) = i(A)\}$. If $H$ is
a $k$-graph or $k$-complex we write $H_A$ for the collection of
sets in $H$ of index $A$. In particular, if $S$ is an $r$-partite
subset of $X$ and $A \sub i(S)$
then $S_A = S \cap \cup_{i \in A} X_i$.
We also use $H_A : K_A(X) \to \{0,1\}$
to denote the characteristic function of this set, i.e.
$H_A(S)$ is $1$ if $S \in H_A$ and $0$ otherwise.

Write $H_A^*$ for the collection of sets $S$ of index $A$ such that all
proper subsets of $S$ belong to $H$. (Note that $H_A \sub H_A^*$.)
We also use $H_A^*$ to denote
the characteristic function of this set.
The relative $A$-density of $H$ is $d_A(H) = |H_A|/|H_A^*|$.
We shorten this to $d_A$ if $H$ is clear from the context.
In particular we have $d_\emptyset=1$,
since $H_\emptyset=H_\emptyset^*=\{\emptyset\}$.
We often assume that
$H_{\{i\}} = X_i$, so that $d_{\{i\}}=1$ for $1 \le i \le r$:
in general we can apply a result obtained under this assumption
by replacing $X = X_1 \cup \cdots \cup X_r$
with $X' = X'_1 \cup \cdots \cup X'_r$,
where $X'_i = H_{\{i\}}$.
We also use an unsubscripted $d$ to denote (absolute) density,
e.g. $d(H_A) = |H_A|/|K_A(X)|$.

Suppose $Y = Y_1 \cup \cdots \cup Y_r$ is another $r$-partite set.
We let $\Phi(Y,X)$ denote the set of all $r$-partite maps
from $Y$ to $X$: these are maps $\phi: Y \to X$ such that
$\phi(Y_i) \sub X_i$ for each $i$. If $J$ is a $k$-graph
or $k$-complex on $Y$ and $H$ a $k$-graph or $k$-complex on $X$
we say that $\phi$ is a homomorphism if $\phi(J) \sub H$. 

For each number $i$ let $U_i = \{u^0_i,u^1_i\}$ be a set of size $2$
and let $U = \cup_i U_i$. For a set of numbers $A$ let
$U_A = \cup_{i \in A} U_i$. The $A$-octahedron on $U$ is
$O_A = \{B \sub U_A: |B \cap U_i|=1, i \in A \}$.
Suppose that for each $B \in O_A$ we have functions
$f_B$ defined on $K_A(X)$. Their Gowers inner product is
$\gen{\{f_B\}_{B \in O_A}}_{\Box^d}
 = \mb{E}_{\omega \in \Phi(U,X)} \prod_{B \in O_A} f_B(\omega(B))$.
(Effectively we are averaging over $\omega$ in $\Phi(U_A,X_A)$,
but extending the function does not affect the average.)
Given a function $f$ defined
on $K_A(X)$ let $\mbox{Oct}(f) = \gen{\{f_B\}}_{\Box^d}$,
where $f_B=f$ for every $B \in O_A$.
Say that $f$ is $\eta$-quasirandom relative to $H$ if
$\mbox{Oct}(f) \le \eta  \prod_{B \in O_A^<} d_B(H)$.\footnote{
Recall that $O_A^<$ is a $(k-1)$-complex consisting of
strict subsets of the sets in $O_A$, and
write $d_B(H) = d_{i(B)}(H)$ for less cumbersome notation.
This will be a general rule: a subscript $A$ should
be understood as $i(A)$ where appropriate.}

If $H$ is a $k$-complex on $X$ and
$A \in \binom{[r]}{\le k} = \{B: B \sub [r], |B| \le k\}$,
the balanced function of $H_A$ is
$\ov{H_A} = (H_A - d_A)H_A^* = H_A - d_AH_A^*$.
We say that $H_A$ is $\eta$-quasirandom if its
balanced function $\ov{H_A}$ is $\eta$-quasirandom.
If $F$ is a fixed $r$-partite $k$-complex
we say that $H$ is $(\eps,F,k)$-quasirandom if
$H_A$ is $\eta_{|A|}$-quasirandom for every $A \in F$,
where the `hidden
parameters' $\eta_2, \cdots, \eta_k$ are defined recursively by
$\eps_k=\eps$, and $\eta_s = \frac{1}{2}(\eps_s \prod_{A \in F, |A| \ge s} d_A)^{2^s}$,
$\epsilon_{s-1} = \eta_s |F|^{-1} \prod_{t=s}^k 2^{-t}$ for $2 \le s \le k$.
This terminology is a convenient way of expressing a sufficient condition
for approximate counting of homomorphisms from $F$ to $H$.
It is sometimes helpful to have a notation including the hidden parameters,
thus we say that $H$ is
$\ov{\eta}$-quasirandom, where $\ov{\eta}=(\eta_1,\cdots,\eta_k)$,
if $H_A$ is $\eta_{|A|}$-quasirandom with respect to $H$
for every $A \in \binom{[r]}{\le k}$
(we can take $\eta_1=0$),
and we say that $H$ is $\ov{d}$-dense,
where $\ov{d}=(d_1,\cdots,d_k)$,
if $H_A$ has relative density at least $d_{|A|}$ with respect to $H$
for every $A \in \binom{[r]}{\le k}$.
We also say that $H$ is $(F,\ov{d})$-dense,
if $H_A$ has relative density at least $d_{|A|}$ with respect to $H$
for every $A \in F$. There is a natural hierarchy
$0 \le \eta_1 \ll d_1 \ll \eta_2 \ll \eps_2 \ll d_2 \ll \cdots \ll \eta_k \ll \eps_k
\ll d_k, 1/|F| \le 1$
that arises in applications,
and assuming that the functions implicit in the $\ll$-notation decay sufficiently quickly,
if $H$ is $\ov{\eta}$-quasirandom and $(F,\ov{d})$-dense
then $H_{A^\le}$ is $(\eps_{|A|},F_{A^\le},|A|)$-quasirandom
for every $A \in \binom{[r]}{\le k}$.

\medskip

\nib{Comments.}
We prefer the term $k$-complex ($k$-chain is used in \cite{G2}), as we
are dealing with simplicial complexes. We are using the letter $d$ for
densities ($\delta$ is used in \cite{G2}) to reserve the letter
$\delta$ for minimum degree or density. We should emphasise that
$d$ with a set subscript indicates relative density with respect to
that set, whereas $d$ with no subscript indicates (absolute) density
(which will be generally well-approximated by a product of relative
densities). For the sake of consistency
we always use set subscripts to indicate some kind of restriction to that set.
The functions defined on $K_A(X)$ may be
regarded as $A$-functions as defined in \cite{G2}. The Gowers inner product
was formalised in \cite{GT} in a slightly restricted context and
generally in \cite{Tao2}, where it was called the cube inner product.
Here we are departing from \cite{G2}, where a counting rather than
averaging convention is used, as we find it more convenient not to have
to keep track of normalising factors.
Although we do not need these facts, it may aid the reader to know that
for $d \ge 2$ the operation $f \mapsto \| f \|_{\Box^d} =
\mbox{Oct}(f)^{2^{-d}}$ defines a norm,
and there is a Gowers-Cauchy-Schwartz
inequality $|\gen{\{f_B\}_{B \in O_A}}_{\Box^d}|
\le \prod_{B \in O_A} \| f_B \|_{\Box^d}$.

\subsection{Quasirandom decomposition and the counting lemma}

As for $3$-graphs, the above notion of quasirandomness
admits a decomposition theorem and a counting lemma.

Suppose as before, that $X=X_1 \cup \cdots \cup X_r$ is
an $r$-partite set. A partition $k$-system $P$
is a collection of partitions
$P_A$ of $K_A(X)$ for every $A \in \binom{[r]}{\le k}$.
For each $A$ there is a natural refinement of $P_A$,
called strong equivalence in \cite{G2}, where
sets $S,S'$ in $K_A(X)$ are strongly equivalent
if $S_B$ and $S'_B$ belong to the same class of $P_B$
for every $B \sub A$.\footnote{
Recall that set subscripts indicate an appropriate restriction:
$S_B = S \cap \cup_{i \in B} X_i$.}
Given $x \in K_{[r]}(X)$
the induced complex $P(x)$ has maximal edges equal to
all sets strongly equivalent to some $x_A = \{x_i: i \in A\}$
with $|A|=k$. The main case of interest is when $H$ is a $k$-graph
or $k$-complex and $P_A = (H_A, K_A(X) \sm H_A)$.

The following decomposition theorem is a variant of Theorem 7.3 in \cite{G2},
with two key differences:
(i) we have an equitable partition of $X$, together with an
exceptional class, and
(ii) we want some flexibility in the choice of the hidden
parameters in quasirandomness, rather than the specific choice
inherent in the definition of $(\eps,F,k)$-quasirandomness.
The proof is very similar to that given in \cite{G2}, so we will
defer a sketch of the necessary modifications that give this version
to appendix A.

\begin{theo} \label{decomp}
Suppose $0 < 1/n \ll d_1 \ll \eta_2 \ll d_2 \ll \eta_3 \ll \cdots \ll d_{k-1}
\ll \eta_k \ll \eps, 1/r, 1/m < 1$,
that $X = X_1 \cup \cdots \cup X_r$ is an $r$-partite set with each $|X_i| \ge n$,
and $P$ is a partition $k$-system on $X$ such that each
such that each $P_A$ partitions $K_A(X)$ into at most $m$ sets.
Then there is a partition $k$-system $Q$ refining $P$ such that
\begin{itemize}
\item $Q_A=P_A$ when $|A|=k$,
\item $|Q_A| \le d_{|A|}^{-1/2}$ when $|A|<k$,
\item $Q_i = Q_{\{i\}}$, $1 \le i \le r$ are of the form
$X_i^{0,1} \cup \cdots \cup X_i^{0,m_{i,0}}
\cup X_i^1 \cup \cdots \cup X_i^{m_i}$,
where $X_i^0 = X_i^{0,1} \cup \cdots \cup X_i^{0,m_{i,0}}$
are exceptional classes of size $|X_i^0| < \eps|X_i|$,
and the $X_i^t$, $1 \le i \le r$, $t \ne 0$ are all of equal size, and
\item if $x=(x_1,\cdots,x_r)$ is
chosen uniformly at random from $K_{[r]}(X)$,
then with probability at least $1-\eps$ the induced complex $Q(x)$ is
$\ov{\eta}$-quasirandom and $\ov{d}$-dense, for some $0<d_k<1$.
\end{itemize}
\end{theo}

Next, we have the following generalised counting lemma, the
adaptation of Theorem 5.1 in \cite{G2} to the setting
of $\ov{\eta}$-quasirandom $\ov{d}$-dense
complexes. (The proof is identical.)

\begin{theo} \label{gencount}
Suppose $0 < 1/n \ll d_1 \ll \eta_2 \ll d_2 \ll \eta_3 \ll \cdots \ll d_{k-1}
\ll \eta_k \ll \eps, d_k, 1/|F| < 1$,
$H$ and $F$ are $r$-partite $k$-complexes on
$X = X_1 \cup \cdots \cup X_r$ and $Y = Y_1 \cup \cdots \cup Y_r$ respectively,
and $H$ is $\ov{\eta}$-quasirandom and $(F,\ov{d})$-dense.
Suppose $F_0$ is a subcomplex of $F$, for each $e \in F$ we have a function
$f_e$ on $K_e(X)$ with $f_e=f_eH_e$, and $f_e=H_e$ for all $e \in F \sm F_0$.
Then
$$\mb{E}_{\phi \in \Phi(Y,X)} \left[ \prod_{e \in F} f_e(\phi(e)) \right]
= \prod_{e \in F \sm F_0} d_e(H)
\cdot \mb{E}_{\phi \in \Phi(Y,X)} \left[ \prod_{e \in F_0} f_e(\phi(e)) \right]
\pm \eps \prod_{e \in F} d_e(H).$$
\end{theo}

As for $3$-graphs, there are two special cases of Theorem \ref{gencount}
that are particularly useful. Firstly,
when $F_0 = \emptyset$ we get a counting
lemma for the partite homomorphism density of $F$ in $H$:
$$d_F(H) = \mb{E}_{\phi \in \Phi(Y,X)} \left[ \prod_{e \in F} H_e(\phi(e)) \right]
= (1 \pm \eps) \prod_{e \in F} d_e(H).$$
The further special case when $F$ is a single edge gives an approximation
of absolute densities in terms of relative densities:
$$d(H_A)=(1 \pm \eps)\prod_{B \sub A} d_B(H).$$
Secondly, we see that quasirandomness implies vertex-uniformity
(the proof is almost identical to that given for $3$-graphs,
so we omit it):

\begin{coro} \label{vertex-uniform}
Suppose $0 < 1/n \ll d_1 \ll \eta_2 \ll d_2 \ll \eta_3 \ll \cdots \ll d_{k-1}
\ll \eta_k \ll \eps, d_k < 1$
and $H$ is an $\ov{\eta}$-quasirandom $\ov{d}$-dense
$k$-partite $k$-complex on $X = X_1 \cup X_2 \cup X_3$,
with $H_i=X_i$, $1 \le i \le k$.
Suppose $W_i \sub X_i$, $1 \le i \le r$ and let $H_{[k]}[W]$ be the restriction
of $H_{[k]}$ to $W = W_1 \cup \cdots \cup W_k$.
Then
$$\frac{|H_{[k]}[W]|}{|X_1|\cdots |X_k|} = \frac{|H_{[k]}|}{|X_1|\cdots |X_k|} \cdot
\left( \frac{|W_1| \cdots |W_k|}{|X_1| \cdots |X_k|} \pm \eps \right).$$
\end{coro}

\subsection{The homomorphism complex: quasirandom counting}

Now we will present our quasirandom version of
the counting lemma in full generality.
Suppose $J$ is a $t$-partite $k$-complex with vertex set
$E = E_1 \cup \cdots \cup E_t$ and $G$ is an
$t$-partite $k$-complex with vertex set
$Y = Y_1 \cup \cdots \cup Y_t$.
We define the homomorphism complex
$J \to G$, which we also denote $G'$ for the sake
of compact notation. The vertex set is $Y' = Y'_1 \cup \cdots
\cup Y'_t$, where $Y'_i$ is the set of all maps
$\phi_i: E_i \to Y_i$. For each $A \sub [t]$, the edges $G'_A$
of index $A$ consist of all $|A|$-tuples $(\phi_i)_{i \in A}$
for which the associated $|A|$-partite map $\phi_A: E_A \to Y_A$
is a homomorphism, i.e. $\phi_A(J_A) \sub G_A$.
Note that $G' = J \to G$ is formally a $t$-complex, but
its sets of size bigger than $k$ have a trivial structure,
in that they are present exactly when all their subsets
of size (at most) $k$ are present.

\begin{theo} \label{aux}
Suppose $0 \le 1/n \ll \eta_2 \ll \eps_2 \ll \eta'_2 \ll d_2 \ll \eta_3 \ll \eps_3
\ll \eta'_3 \ll \cdots \ll \eta_k \ll \eps_k \ll \eta'_k \ll 1/t,d_k < 1$,
that $J$ is a $t$-partite $k$-complex on
$E = E_1 \cup \cdots \cup E_t$
and $G$ is an $\ov{\eta}$-quasirandom $(J,\ov{d})$-dense
$t$-partite $k$-complex
on $Y = Y_1 \cup \cdots \cup Y_t$, where $G_i=Y_i$, $1 \le i \le t$.
Then $G' = J \to G$ is $\ov{\eta'}$-quasirandom,
where $\ov{\eta'}=(\eta'_2,\cdots,\eta'_t)$ with $\eta'_i=0$ for $k < i \le t$.
Furthermore, the densities may be
estimated as
$d_A(G') = (1 \pm \eps) d_A(G)^{|J_A|}$
for $|A| \le k$ and $d_A(G') = 1$ for $|A|>k$.
\end{theo}

\nib{Proof.}
The argument is very similar to that used in the proof of Theorem \ref{qr-oct}.
We write $d_A=d_A(G)$ and $d'_A=d_A(G')$ for $A \sub [t]$, $|A| \ge 2$.
Fix $A \sub [t]$, $|A| \ge 2$. The counting lemma gives
$$d'_A  =   \frac{|G'_A|}{|G^{'*}_A|}
 = \frac{\mb{P}_{\phi \in \Phi(E,Y)}(\phi(J_A^\le) \sub G)}{
  \mb{P}_{\phi \in \Phi(E,Y)} (\phi(J_A^<) \sub G)}
 = \frac{(1 \pm \eps_{|A|}/3) \prod_{S \in J_A^\le} d_S}{
   (1 \pm \eps_{|A|}/3) \prod_{S \in J_A^<}  d_S}
 = (1 \pm \eps_{|A|}) d_A^{|J_A|}.$$
Next we need to estimate $\mbox{Oct}(\ov{G'_A})$.
Write $\ov{G'_A}=g'_0-g'_1$,
where $g'_0 = G'_A$ and $g'_1 = d'_A G^{'*}_A$.
We can assume that $|A| \le k$:
otherwise we have $G'_A = G^{'*}_A$, $d'_A=1$
and $\ov{G'_A} = 0$.
Then
$$\mbox{Oct}(\ov{G'_A}) =
\mb{E}_{\omega \in \Phi(U,Y')} \prod_{B \in O_A} \ov{G'_A}(\omega(B))
= \sum_j (-1)^{\sum j}\ \mb{E}_{\omega \in \Phi(U,Y')}
\prod_{B \in O_A} g'_{j_B}(\omega(B)),$$
where we recall that $O_A$ is the $A$-octahedron on $U$,
$j$ ranges over all $\{0,1\}$-sequences $(j_B)_{B \in O_A}$,
and $\sum j = \sum_{B \in O_A} j_B$.

We interpret each summand as counting homomorphisms to functions from $JJ_A$,
where we define the clone $JJ$ of $J$ as the following $t$-partite $k$-complex
on $EE = EE_1 \cup \cdots \cup EE_t$. The parts
$EE_i = E_i^0 \cup E_i^1$, $1 \le i \le t$ have two copies $x_i^0, x_i^1$
of each point $x_i \in E_i$. The edges of $JJ_A$ consist of $2^{|A|}$ copies of $J_A$,
indexed by $O_A$, thus $JJ_A = \cup_{B \in O_A} J_A^B$,
writing $J_A^B = \{ \{x_i^B\}_{i \in A}: \{x_i\}_{i \in A} \in J_A \}$,
where $x_i^B$ is $x_i^0$ if $B \cap U_i = \{u_i^0\}$ or $x_i^1$ if
$B \cap U_i = \{u_i^1\}$.
Given a $t$-partite map $\omega$ from $U$ to $Y'$,
we obtain $t$-partite maps $\phi^B=\omega(B)$, $B \in O_A$
acting on various copies of $E$ in $EE$,
which together give a $t$-partite map $\phi$ from $EE$ to $Y$
defined by $\phi(x_i^j)=\omega(u_i^j)(x_i)$.
Now we can rewrite the summands above using
$$\prod_{B \in O_A} g'_{j_B}(\omega(B))
= \prod_{S \in O_A} \prod_{T \in J^B} g_{j_B}(\phi^B(T))
= \prod_{T \in JJ_A} g_{j_{B(T)}}(\phi(T)),$$
where we write $B(T)$ for that $B \in O_A$ such that $Y \in J_A^B$,
$g_0 = G_A$ (the characteristic function) and
$g_1 = (d'_A)^{1/|J_A|}G^*_A = (1 \pm 2|J_A|^{-1}\eps_{|A|})d_A G^*_A$.

Now we can estimate the contribution from the summand corresponding to
$j=(j_B)_{B \in O_A}$ using the counting lemma for
homomorphisms to $G$ from the complex
$K_j = JJ_{A^\le} \sm \cup_{B:j_B=1} J_A^B$.
We have
\begin{align*}
\mb{E}_{\omega \in \Phi(U,Y')} \prod_{B \in O_A} g'_{j_B}(\omega(B))
& = \mb{E}_{\phi \in \Phi(EE,Y)} \prod_{T \in JJ_A} g_{j_{B(T)}}(\phi(T)) \\
& = (d'_A)^{\sum j} \mb{P}_{\phi \in \Phi(EE,Y)}[\phi(K_j) \sub G] \\
& = ((1 \pm \eps_{|A|})d_A^{|J_A|})^{\sum j}
\cdot (1 \pm \eps_{|A|}) \prod_{T \in K_j} d_T(G) \\
& = (1 \pm O(\eps_{|A|})) \prod_{T \in JJ_{A^\le}} d_T(G)
 = (1 \pm O(\eps_{|A|})) \prod_{B \in O_{A^\le}} d'_B.
\end{align*}
Therefore
$$\mbox{Oct}(\ov{G'_A}) =
= \sum_j (-1)^{\sum j} (1 \pm O(\eps_{|A|})) \prod_{B \in O_{A^\le}} d'_B
= \sum_j  O(\eps_{|A|}) \prod_{B \in O_{A^\le}} d'_B
< \eta'_{|A|} \prod_{B \in O_{A^<}} d'_B,$$
since $\eps_{|A|} \ll \eta'_{|A|}$.
This proves that $G'_A$ is $\eta'_{|A|}$-quasirandom, so we are done. \qed

\subsection{An application to generalised Tur\'an problems}

In this subsection we apply the quasirandom counting lemma to derive
information about generalised Tur\'an problems for configurations
that have the following particular structure.
Suppose $F$ is a $k$-graph and $s$ is an integer. The $s$-blowup
$F(s)$ of $F$ is defined as follows. For each vertex $x$ of $F$ there
are vertices $x_1,\cdots,x_s$ of $F(s)$. For each edge $x^1 \cdots
x^k$ of $F$ we have all $s^k$ edges $x^1_{i_1} \cdots x^k_{i_k}$ with
$1 \le i_1,\cdots,i_k \le s$ in $F(s)$. Note that $F(2)=FF$ is the clone
of $F$ as defined in the previous section. For $1 \le s \le k-1$ we define
the $s$-augmentation $F^{+s}$ of $F$ as the $k$-graph obtained from
$F(s)$ by adding a set $V^+$ of $k-s$ new vertices and all edges
$V^+ \cup \{x_1,\cdots,x_s\}$, with $x \in V(F)$.
The following theorem is a generalisation of Theorem \ref{sparse3} in two respects:
(i) the relationship between a single edge and the Fano plane is replaced
by the relationship between $F$ and $F^{+s}$,
(ii) the assumptions have been relaxed to `$\theta$-approximate' assumptions.

\begin{theo} \label{structure}
Suppose $0 < 1/n \ll \theta \ll \eps \ll 1/k$, $1 \le s \le k-1$, $0 < \delta < 1$,
$F$ is a $k$-graph, $H$ is a $k$-graph on a set $X$ of $n$ vertices,
$d_{F^{+s}}(H) < \theta$ and the $s$-graph
$G = \{S \in \binom{X}{s}: |N_H(S)| <  \delta \binom{n-s}{k-s} \}$
has density $d(G) < \theta$.
Then there is a subset $Z \sub X$ with $|Z| > (\delta-\eps)n$ so that
$d_F(H[Z]) < \eps$.
\end{theo}

\nib{Proof of Theorem \ref{structure}.}
We operate with a constant hierarchy
$0 < 1/n \ll \theta \ll \alpha
\ll d_1 \ll \eta_2 \ll d_2 \ll \eta_3 \ll \cdots \ll d_{k-1} \ll \eta_k \ll \gamma, d_k
\ll \omega \ll \eps \ll 1$.
Set $r=\lceil \omega^{-2} \rceil$ and let $X = X_1 \cup \cdots \cup X_r$ be an arbitrary partition
with $||X_i| - n/r| < 1$ for $1 \le i \le r$.
Let $H_1$ be the edges of $H$ that respect this partition, i.e.
have at most one point in each $X_i$.
Consider the partition system $P$ naturally associated with $H_1^\le$, in which
$P_A$ is the partition $K_A(X) = (H_1^\le)_A \cup K_A(X) \sm (H_1^\le)_A$.
Applying Theorem \ref{decomp}, using $\gamma$ instead of $\eps$,
we obtain a partition system $Q$ refining $P$ such that
\begin{itemize}
\item $Q_A=P_A$ when $|A|=k$,
\item $|Q_A| \le d_{|A|}^{-1/2}$ when $|A|<k$,
\item $Q_i = Q_{\{i\}}$, $1 \le i \le r$ are of the form
$X_i^{0,1} \cup \cdots \cup X_i^{0,m_{i,0}}
\cup X_i^1 \cup \cdots \cup X_i^{m_i}$,
where $X_i^0 = X_i^{0,1} \cup \cdots \cup X_i^{0,m_{i,0}}$
are exceptional classes of size $|X_i^0| < \gamma|X_i|$,
and the $X_i^t$, $1 \le i \le r$, $t \ne 0$ are all of equal size, say $m$, and
\item if $x=(x_1,\cdots,x_r)$ is
chosen uniformly at random from $K_{[r]}(X)$,
then with probability at least $1-\gamma$ the induced complex $Q(x)$ is
$\ov{\eta}$-quasirandom and $\ov{d}$-dense, for some $0<d_k<1$.
\end{itemize}

Consider the set of all $s$-tuples contained in the vertex classes of $Q$:
denote it $Y = \cup_{(i,t_i) \in I} Y_{i,t_i}$,
where $I=\{(i,t_i): 1 \le i \le r, 1 \le t_i \le n_i\}$
and $Y_{i,t_i} = \binom{X_i^{t_i}}{s}$.
Let $W = \binom{X}{k-s}$ and
$W' = \{R \in W: |N_H(R) \cap Y| > (\delta-\omega)|Y| \}$.
Double counting gives
$$(|Y|-\theta n^s/s!) \delta \binom{n-s}{k-s} <
\sum_{S \in Y \sm G} d_H(S)
 = \sum_{R \in W} |N_H(R) \cap Y \sm G| <
|W'||Y| + |W \sm W'|(\delta-\omega)|Y|,$$
which gives $|W'| > \frac{1}{2} \omega \binom{n}{k-s}$.
Write
$C_R = \{(i,t_i) \in I:  |N_H(R) \cap Y_{i,t_i}| > \omega \binom{m}{s} \}$.
Since $\alpha \ll d_1$
there are at most $2^{|I|} < \alpha^{-1/2}$ possibilities for $C_R$,
so there is some $C$ that occurs as $C_R$ for
at least $\alpha \binom{n}{k-s}$
of the $R$ in $W'$. For now we focus on one such $R$.
Then, since $|X_i^{t_i}|=m$ for every $(i,t_i) \in I$,
$$(\delta-\omega)|I|\binom{m}{s} < |N_H(R) \cap Y| <
\sum_{(i,t_i) \in C} \binom{m}{s} + \sum_{(i,t_i) \in I \sm C}
\omega \binom{m}{s},$$
so $|C| > (\delta-2\omega)|I|$. Let $Z = \cup_{(i,t_i) \in C} X_i^{t_i}$.
Then $|Z| = |C|m > (\delta-3\omega)n$.

Suppose $V(F) = \{v^1,\cdots,v^t\}$ for some $t$.
We think of $F^\le$ as a $t$-partite $k$-complex with one vertex in each part.
For any (ordered) $t$-tuple of vertices $x=(x_1,\cdots,x_t)$
with $x_a \in X_{i_a}^{t_a}$, $(i_a,t_a) \in C$, $i_a$ distinct, $1 \le a \le t$
we consider the induced $t$-partite $k$-complex $Q(x)$,
i.e. the union of the induced complexes corresponding to each $k$-subset of $x$.
For convenient notation we temporarily identify $(i_1,\cdots,i_t)$ with $[t]=\{1,\cdots,t\}$.
We let $H(x) \sub Q(x)$ be obtained by keeping those $k$-tuples corresponding to edges of $H$,
i.e. $H(x)_A = Q(x)_A$ if $Q(x)_A \sub H_A$ and $H(x)_A = \emptyset$ if
$Q(x)_A \sub K_A \sm H_A$, for $A \in \binom{[t]}{k}$.
We claim that $H(x)$ is not both $\ov{\eta}$-quasirandom and $(F^\le,\ov{d})$-dense.
Otherwise, we can calculate as follows that there is too high a probability
that a random map $\phi:V(F^{+s})\to X$ is a homomorphism from $F^{+s}$ to $H$.
Using the definition of $F^{+s}$,
we write the event $\{\phi(F^{+s}) \sub H\}$ as $E_1 \cap E_2 \cap E_3 \cap E_4$,
where $E_1$ is the event that the additional vertex set $V^+$
is mapped to some $R$ in $W'$ with $C_R=C$,
$E_2$ is the event that $V(F(s))$ is mapped in $r$-partite fashion
to the classes  $X_{i_a}^{t_a}$, $1 \le a \le t$,
$E_3$ is the event that each $s$-tuple $v^a_1 \cdots v^a_s$
is mapped to $N_H(R) \cap Y_{i_a,t_a}$ for $1 \le a \le t$,
and $E_4$ is the event $\{\phi(F(s)) \sub H\}$.
Now $\mb{P}(E_1)>\alpha/2$ by definition of $C$,
$\mb{P}(E_2|E_1)=(m/n)^{ts}>\alpha$ and
$\mb{P}(E_3|E_1 \cap E_2) > (\omega/2)^t$ by definition of $C_R$.
Also, Theorem \ref{aux} and Corollary \ref{vertex-uniform}
applied to $J=F(s)^\le$ give
$$\mb{P}(E_4|E_1 \cap E_2 \cap E_3)
= (1 \pm \eps) \prod_{A \in F(s)^\le} d_A(H(x))
> \frac{1}{2} \prod_{A \in F^\le} d_{|A|}^{2^{|A|}} > \alpha.$$
We deduce that $\mb{P}[\phi(F^{+s}) \sub H] > \alpha/2 \cdot \alpha
\cdot (\omega/2)^t \cdot \alpha > \theta$, which contradicts
our assumptions. Therefore $H(x)$ is not both $\ov{\eta}$-quasirandom
and $(F^\le,\ov{d})$-dense.

Now we want to estimate $d_F(H[Z]) = \mb{P}_{f:E \to Z}(f(F) \sub H[Z])$.
The event $f(F) \sub H[Z]$ is contained in the event
$B_1 \cup B_2 \cup B_3$ where $B_1$ is the event that $f$ is not a partite map
(i.e. two vertices of $F$ are mapped to the same $X_i$),
$B_2$ is the event that $f$ is a partite map but $F$ is mapped to $H(x)$
for some $t$-tuple $x$ that is not $\ov{\eta}$-quasirandom,
and $B_3$ is the event that neither $B_1$ or $B_2$ holds and $f(F) \sub H$.
Then $\mb{P}(B_1) < \binom{t}{2}/r < \omega$.
Also, $Q$ was chosen so that at most $\gamma n^k$
edges of $H_1$ belong to induced complexes which fail a quasirandomness condition,
so $\mb{P}(B_2) < |E(F)| \gamma < \omega$.
To estimate $\mb{P}(B_3)$ we note that on this event
$f$ maps $V(F)$ to some $H(x)$ that is $\ov{\eta}$-quasirandom, so by the previous paragraph
it is not $(F^\le,\ov{d})$-dense. Since $|Q_A| \le d_{|A|}^{-1/2}$ when $|A|<k$, for any $a<k$,
at most $\binom{k}{a}d_a^{1/2}n^k$ edges of $H_1$ belong to an induced complex
in which some $a$-graph has relative density at most $d_a$.
Also, at most $d_k n^k$ edges of $H_1$ belong to an induced complex in which
the $k$-graph has relative density at most $d_k$. Summing these contributions,
we deduce that $\mb{P}(B_3) < \omega$.
Finally we have $d_F(H[Z]) < \mb{P}(B_1) + \mb{P}(B_2) + \mb{P}(B_3) < 3\omega < \eps$,
as required. \qed

\medskip

\nib{Remark.}\
As a further extension, note that we could have
taken any $I' \sub I$ and applied the argument with
$Y' = \cup_{(i,t_i) \in I'} Y_{i,t_i}$ instead of $Y$. Thus we find some
$C' \sub I'$ with $|C'|>(\delta-\eps)|I'|$ such that
$Z' = \cup_{(i,t_i) \in C'} X_i^{t_i}$
satisfies $|Z'|=|C'|m > \frac{|I'|}{|I|} (\delta-\eps)n$
and $d_F(H[Z']) < \eps$.

\section{Codegree problems for projective planes}

In this section we apply Theorem \ref{structure} to codegree problems
for projective planes.
The first subsection contains definitions and a summary of previous results
on codegree problems for projective planes obtained in \cite{KZ}.
The second subsection generalises the approach used for the Fano plane to obtain
strong structural information for the general problem,
which already determines the codegree threshold for planes over
a field of odd size up to an additive constant
(the first part of Theorem \ref{pg2-odd}).
In the third subsection we complete the proof of Theorem \ref{pg2-odd}
with a more detailed analysis for the plane over $\mb{F}_3$ in terms
of its blocking sets: we see that the true codegree threshold in this case
lies strictly between the natural upper and lower bounds found in the first subsection.
The last subsection deals with the plane over $\mb{F}_4$, where we demonstrate
a surprisingly different behaviour from $\mb{F}_2$ and fields of odd size.

\subsection{Definitions and previous results}

First we give some definitions.
Let $\mb{F}_q$ be the field with $q$ elements, for any
prime power $q$.
The projective geometry $PG_m(q)$ of dimension $m$ over $\mb{F}_q$
is a $(q+1)$-graph with vertex set equal to the one-dimensional subspaces of
$\mb{F}_q^{m+1}$ and edges corresponding to the two-dimensional subspaces
of $\mb{F}_q^{m+1}$, in that for each
two-dimensional subspace, the set of one-dimensional
subspaces that it contains is an edge of the hypergraph $PG_m(q)$.
A vertex of $PG_m(q)$ can be described by projective
co-ordinates as $(x_1:\cdots:x_{m+1})$, where $(x_1,\cdots,x_{m+1})$ is
any non-zero vector of $\mb{F}_q^{m+1}$ and
$(x_1:\cdots:x_{m+1})$ denotes the one-dimensional subspace
that it generates.

The main result of \cite{KZ} is the following
upper bound on the codegree density for general projective
geometries, which is tight in many cases.

\begin{theo}
The codegree density of projective geometries satisfies
$\pi_q(PG_m(q)) \le 1-1/m$. Equality holds whenever $m=2$ and $q$ is $2$ or odd,
and whenever $m=3$ and $q$ is $2$ or $3$.
\end{theo}

The results in \cite{KZ} can be summarised by the following table, in which the entry
in the cell indexed by row $m$ and column $q$ is either a number
indicating the exact value of $\pi_q(PG_m(q))$ or an interval
in which $\pi_q(PG_m(q))$ lies.

\medskip

\begin{tabular}{|c|ccccc|}
\hline
$m \backslash q$ & $2$ & $3$ & $4$ & $2^t$, $t \ge 3$ & $p^t$, $p$ odd \\
\hline
$2$ & $1/2$ & $1/2$ & $[1/3,1/2]$ & $[0,1/2]$ & $1/2$ \\
$3$ & $2/3$ & $2/3$ & $[1/2,2/3]$ & $[0,2/3]$ & $[1/2,2/3]$ \\
$4$ & $[2/3,3/4]$ & $[2/3,3/4]$ & $[1/2,3/4]$ & $[0,3/4]$ & $[1/2,3/4]$ \\
$m \ge 5$ & $[3/4,1-1/m]$ & $[2/3,1-1/m]$ & $[1/2,1-1/m]$ & $[0,1-1/m]$ &
$[1/2,1-1/m]$ \\
\hline
\end{tabular}

\subsection{General structure}

A useful property of projective geometries proved in \cite{KZ} is that
$PG_m(q) \sub PG_{m-1}(q)^{+q}$; in particular
$PG_2(q) \sub e^{+q}$, where $e$ is a single $(q+1)$-edge.
Thus we may specialise Theorem \ref{structure} as follows.

\begin{theo} \label{pg2-structure}
Suppose $q$ is a prime power, $0 < 1/n \ll \theta \ll \eps \ll 1/q$, $0<\delta<1$,
$H$ is a $(q+1)$-graph on a set $X$ of $n$ vertices,
$d_{PG_2(q)}(H) < \theta$ and the $q$-graph
$G = \{Q \in \binom{X}{q}: |N_H(Q)| < \delta n \}$
has density $d(G) < \theta$.
Then there is a subset $Z \sub X$ with $|Z| > (\delta-\eps)n$ so that
$d(H[Z]) < \eps$.
\end{theo}

Similarly to Theorem \ref{stab3},
we can use this to deduce a stability result, describing the
approximate structure of a $(q+1)$-graph that
does not contain $PG_2(q)$ and has most of its
$q$-degrees at least $n/2-o(n)$.
\footnote{This weakened assumption may
strike the reader as strange at first sight, as in a straightforward
application of the stability method one has information
about all $q$-degrees, but its importance will become clearer
in the third subsection.}

\begin{theo} \label{pg2-stab}
Suppose $q$ is a prime power, $0 < 1/n \ll \theta \ll \eps \ll 1$,
$H$ is a $(q+1)$-graph on a set $X$ of $n$ vertices,
$d_{PG_2(q)}(H) < \theta$ and the $q$-graph
$G = \{Q \in \binom{X}{q}: |N_H(Q)| < (1/2-\theta)n \}$
has density $d(G) < \theta$.
Then there is a partition $X = A \cup B$
such that at most $\eps n^{q+1}$ edges are contained
entirely within $A$ or within $B$.
\end{theo}

Before giving the proof we remark that any $k$-graph $H$
with density $d(H) > 1 - \binom{m}{k}^{-1}$ contains
$K^k_m$, the complete $k$-graph on $m$ vertices.
Indeed, by averaging $H$ contains a set $M$ of $m$ vertices
with $d(H[M]) \ge d(H)$, and then every $k$-tuple in $M$
must be an edge. (In other words we are using the
easy bound $\pi(K^k_m) \le 1 - \binom{m}{k}^{-1}$
on the Tur\'an density, which is far from being best possible,
but suffices for our purposes here.)

\medskip

\nib{Proof.} Introduce a hierarchy of constants
$\eps = \theta_4 \gg \theta_3 \gg \theta_2 \gg \theta_1 = \theta$.
By Theorem \ref{pg2-structure} we can find $A \sub X$ with
$|A| > (1/2-\theta_2)n$ and $d(H[A]) < \theta_2$, i.e.
$|H[A]| < \theta_2 |A|^{q+1}/(q+1)!$.
Let $B = X \sm A$. Now
$$(q+1)|H[A]| + \sum_{Q \in \binom{A}{q}} |N_H(Q) \cap B|
= \sum_{Q \in \binom{A}{q}} |N_H(Q)|
> \left| \binom{A}{q} \sm G \right| (1/2-\theta)n
> (1/2 - \theta_2)n\binom{|A|}{q}$$
so
$$\sum_{b \in B} \left| N_H(b) \cap \binom{A}{q}\right|
= \sum_{Q \in \binom{A}{q}} |N_H(Q) \cap B| > (1/2-\theta_2)n \binom{|A|}{q}
- \theta_2 |A|^{q+1}/q! > (1/2 - \theta_3)n \binom{|A|}{q}.$$
This implies that $|B| > (1/2-\theta_3)n$, so
$|A|$ and $|B|$ are $(1/2 \pm \theta_3)n$.
Let
$$B_0 = \left\{b \in B: \left|N_H(b) \cap \binom{A}{q}\right| <
(1-\theta_3^{1/2}) \binom{|A|}{q} \right\}.$$
Then
\begin{eqnarray*}
(1/2 - \theta_3)n \binom{|A|}{q} & < &
\sum_{b \in B} \left| N_H(b) \cap \binom{A}{q} \right|
< |B_0| (1-\theta_3^{1/2}) \binom{|A|}{q} + |B \sm B_0| \binom{|A|}{q} \\
& = & |B| \binom{|A|}{q} - |B_0| \cdot \theta_3^{1/2} \binom{|A|}{q},
\end{eqnarray*}
so
$$|B_0| < \theta_3^{-1/2}(|B| - (1/2 - \theta_3)n) < 2\theta_3^{1/2}n.$$

We claim that there is no edge $e \in H[B \sm B_0]$.
If there were we would have $|\bigcap_{b \in e} N_H(b) \cap \binom{A}{q}|
> (1 - (q+1)\theta_3^{1/2}) \binom{|A|}{q}$. But then, by the remark
before the proof
$\bigcap_{b \in e} N_H(b) \cap \binom{A}{q}$ contains a $K^q_{q^2}$,
and together with $e$ we have a copy of $PG_2(q)$.
Therefore every edge of $B$ contains a point of $B_0$, so
$|H[B]| \le |B_0| \binom{|B|}{q} < 2\theta_3^{1/2}n
((1/2+\theta_3)n)^q/q! < \frac{1}{2} \eps n^{q+1}$.
Since $|H[A]| < \theta_2 |A|^{q+1}/(q+1)!$ we are done. \qed

\medskip

Finally, just as in the proof of Theorem \ref{fano},
we use the previous stability result to obtain very
precise information about a
$(q+1)$-graph on a set $X$ of $n>n_0$ vertices
with minimum $q$-degree $\delta_q(H) \ge n/2$
that does not contain $PG_2(q)$.

\begin{theo} \label{pg2-pre-exact}
For any prime power $q$ there is a number $n_0$ so that if $H$ is a
$(q+1)$-graph on a set $X$ of $n>n_0$ vertices
with minimum $q$-degree $\delta_q(H) \ge n/2$
that does not contain $PG_2(q)$
then $n$ is even, there is a partition $X = X_0 \cup X_1$ where
$|X_0|=|X_1|=n/2$ and $H[X_0]=H[X_1]=\emptyset$.
\end{theo}

\nib{Proof.}
Let $\eps$ be sufficiently small and $n > n_0(\eps,q)$ sufficiently large.
By Theorem \ref{pg2-stab} we have a partition $X = X_0 \cup X_1$
so that  $|H[X_0]|+|H[X_1]| < \eps n^{q+1}$. Choose the partition
that minimises $|H[X_0]|+|H[X_1]|$. We will show that
this partition satisfies the conclusion of the theorem.
Note first that the same argument used in the proof of Theorem \ref{pg2-stab}
shows that $|X_0|$ and $|X_1|$ are $(1/2 \pm 2\eps)n$.

First we show that there is no vertex $x \in X_0$ with degree at least
$\eps^{1/2q} n^q$ in $H[X_0]$. For suppose there is such a vertex $x$.
By choice of partition we have $|N_H(x) \cap \binom{X_1}{q}|
\ge|N_H(x) \cap \binom{X_0}{q}| \ge \eps^{1/2q}n^q$, or we could
reduce $|H[X_0]|+|H[X_1]|$ by moving $x$ to $X_1$.
We can choose matchings
$M = \{x^1_1 \cdots x^q_1,\cdots, x^1_m \cdots x^q_m\}$
in $N_H(x) \cap \binom{X_0}{q}$
and $M' = \{y^1_1 \cdots y^q_1,\cdots,y^1_m \cdots y^q_m\}$
in $N_H(x) \cap \binom{X_1}{q}$,
with $m = \eps^{1/2q} n/q$
(as in the case $q=2$ we are using the well-known observation that any maximal matchings
will be at least this large).
Now
$$ \sum_{I = \{i_1,\cdots,i_q\} \in \binom{[m]}{q}}
\sum_{J = \{j_1,\cdots,j_q\} \in [q]^q}
 |N_{H[X_0]}(x^{j_1}_{i_1} \cdots x^{j_q}_{i_q})|
< (q+1)|H[X_0]| < (q+1)\eps n^{q+1},$$
so we can choose $I \in \binom{[m]}{q}$ such that
$$ \sum_{J = \{j_1,\cdots,j_q\} \in [q]^q}
|N_{H[X_0]}(x^{j_1}_{i_1} \cdots x^{j_q}_{i_q})|
< \binom{m}{q}^{-1} (q+1)\eps n^{q+1} < 2q^q(q+1)! \eps^{1/2} n.$$

Also, for every $1 \le k  \le m$, $H$ cannot have
all $q^{q+1}$ edges $x^{j_1}_{i_1} \cdots x^{j_q}_{i_q} y_k^{j_{q+1}}$, with
$j_1,\cdots,j_{q+1}$ in $[q]$, as then together with $x$
we have a $q$-augmented edge,
which contains $PG_2(q)$. Therefore there is
some $q$-tuple $Q = x^{j_1}_{i_1} \cdots x^{j_q}_{i_q}$ such that
there are at least $q^{-q}m$ vertices $y_k^{j_{q+1}}$ that
do not belong to $N_H(Q)$.
This gives
$|N_H(Q)| < 2q^q(q+1)! \eps^{1/2} n + |X_1|
- q^{-q-1} \eps^{1/2q} n < (1/2-\eps)n$,
which contradicts our assumptions. We deduce that there
is no vertex in $X_0$ with degree at least $\eps^{1/2q} n^q$ in $H[X_0]$.
Similarly there is no vertex in $X_1$
with degree at least $\eps^{1/2q} n^q$ in $H[X_1]$.

Write $|X_0|=n/2+t$ and $|X_1|=n/2-t$, where without loss of
generality $0 \le t \le 2\eps n$. Suppose for a contradiction that either
$t>0$ or $t=0$ and there is an edge in $H[X_0]$ or $H[X_1]$.
Note that for every $q$-tuple $Q \in \binom{X_0}{q}$
we have $|N_H(Q) \cap X_0| \ge |N_H(Q)| - |X_1| \ge t$.
Thus we can assume there is at least one edge in $X_0$ (since the case
$t=0$ is symmetrical). Let $Y_i \sub X_i$ be minimum size
transversals of $H[X_i]$. Then $Y_0 \ne \emptyset$.
Also, by the previous paragraph
$$|Y_0|\eps^{1/2q}n^q > \sum_{x \in Y_0} |N_{H[X_0]}(x)| \ge |H[X_0]|
= \frac{1}{q+1} \sum_{Q \in \binom{X_0}{q} } |N_{H[X_0]}(Q)|
\ge \frac{t}{q+1} \binom{|X_0|}{q}$$
so $|Y_0| > \eps^{-1/2q} \frac{t}{q+1} n^{-q} \binom{n/2+t}{q}
> 2q^{3q} t$, say.

For each edge $e \in X_0$ consider all possible ways to extend it
to a copy of $PG_2(q)$ using some $F \in \binom{X_1}{q^2}$.
Since $H$ does not contain $PG_2(q)$ there is some $(q+1)$-tuple
with $q$ points in $F$ and one point in $e$ that is not an edge.
We count each such $(q+1)$-tuple at most $\binom{|X_1|-q}{q^2-q}$ times, so
we get a set of at least $\binom{|X_1|-q}{q^2-q}^{-1} \binom{|X_1|}{q^2}$
distinct $(q+1)$-tuples.
Thus there is some point $x$ in $e$ for which we have
a set $M_x$ of at least
$\frac{1}{q+1} \binom{|X_1|-q}{q^2-q}^{-1} \binom{|X_1|}{q^2}
> q^{-3q} \binom{|X_1|}{q}$ `missing' $(q+1)$-tuples involving
$x$ and a $q$-tuple in $\binom{X_1}{q}$ not
belonging to $N_H(x)$. Now varying
$e$ over all edges in $X_0$ we get a set $M = \cup_x M_x$ of missing
$(q+1)$-tuples with
$|M| \ge |Y_0| \cdot q^{-3q} \binom{|X_1|}{q}$ (since $Y_0$ is
a minimum size transversal).
For each $q$-tuple $Q \in \binom{X_1}{q}$ we have
$$n/2 \le |N_H(Q)|  = |N_H(Q) \cap X_0| + |N_H(Q) \cap X_1|
\le |X_0|-|N_M(Q)| + |N_H(Q) \cap X_1|$$ so
$$(q+1)|H[X_1]| = \sum_{Q \in \binom{X_1}{q}} |N_H(Q) \cap X_1|
 \ge \sum_{Q \in \binom{X_1}{q}} (|N_M(Q)|-t)
 = |M| - t\binom{|X_1|}{q}.$$
Since $q^{-3q}|Y_0| - t > \frac{1}{2}q^{-3q}|Y_0|$ we have
$|H[X_1]| > \frac{q^{-3q}}{2(q+1)} |Y_0|\binom{|X_1|}{q}$.
Also $|Y_1|\eps^{1/2q}n^q > \sum_{x \in Y_1} |N_{H[Y_1]}(x)| > |H[X_1]|$,
so $|Y_1| >  \eps^{-1/2q}  \frac{q^{-3q}}{2(q+1)} |Y_0| n^{-q} \binom{n/2-t}{q}
> 2|Y_0|$.

Finally we can apply the argument of the previous paragraph
interchanging $X_0$ and $X_1$. We get a set $M'$ of at
least $|Y_1| \cdot q^{-3q} \binom{|X_0|}{q}$ distinct $(q+1)$-tuples
that are not edges, each having $q$ points in $X_0$ and
$1$ point in $X_1$. For each $q$-tuple $Q \sub X_0$ we now have
$|N_H(Q) \cap X_0| \ge |N_{M'}(Q)| + t$, so
$|H[X_0]| > \frac{1}{q+1} \cdot (q^{-3q}|Y_1|+t)\binom{|X_0|}{q}$
and
$|Y_0| > \eps^{-1/2q}n^q |H[X_0]| > 2|Y_1|$.
This contradiction completes the proof. \qed

\medskip

\nib{Remark.} The argument applies more generally to any $F$ with
$PG_2(q) \sub F \sub e^{+q}$ (where $e=PG_1(q)$ is a single edge.)

\medskip

For general odd $q$ Theorem \ref{pg2-pre-exact} determines
the $q$-degree threshold to find $PG_2(q)$ up to a constant,
thus proving the first part of Theorem \ref{pg2-odd}. It provides an upper bound
$\mbox{ex}_q(n,PG_2(q)) \le n/2$. On the other hand, it was
proved in \cite{KZ} that there is no copy of $PG_2(q)$ in the complete
oddly bipartite $(q+1)$-graph, by which we mean the construction
obtained by forming a balanced partition $X=X_0 \cup X_1$
and taking as edges all $(q+1)$-tuples with an odd number of points
in each $X_i$. This gives a lower bound
$\mbox{ex}_q(n,PG_2(q)) \ge \lfloor n/2 \rfloor - q + 1$.

\subsection{The projective plane over $\mb{F}_3$}

To nail down the constant uncertainty in the bound
$\mbox{ex}_q(n,PG_2(q)) \ge \lfloor n/2 \rfloor - q + 1$ for odd $q$
requires analysis of the degrees of $q$-tuples not contained
in $X_0$ or $X_1$, which is closely connected to the theory
of blocking sets in projective planes (see \cite{Tal}).
This theory is far from complete, but the case $q=3$
is sufficiently simple to analyse.

Say that $S \subset PG_m(q)$ is a blocking set if $0<|S \cap L| < |L|$ for
every line $L$ of $PG_m(q)$. Note that the complement of a blocking set
is also a blocking set, so the existence of a blocking
set is equivalent to the existence of a bipartition of $PG_m(q)$.
The blocking sets of $PG_2(3)$ may be classified as follows
(see \cite{Tal} or \cite{H}):
they all have size $6$ or $7$, and those of size $6$ are
of the form $L(x,y) \cup L(y,z) \cup L(x,z) \sm \{x,y,z\}$,
where $x,y,z$ are non-collinear points and we use the notation
$L(a,b)$ for the line containing $a$ and $b$. (Those of size $7$ are their
complements.)

\noindent \parbox{0.5\textwidth}{
Consider a partition $V(PG_2(3)) = A_0 \cup A_1$
of $PG_2(3)$ into blocking sets, where
$A_1 = L(x,y) \cup L(y,z) \cup L(x,z) \sm \{x,y,z\}$.
We refer to the type of a line $L$ as $A_0^{t_0} A_1^{t_1}$ if
$|L \cap A_i|=t_i$, $i=0,1$. Then
the lines $L(x,y)$, $L(y,z)$ and $L(x,z)$ each have
type $A_0^2 A_1^2$, and their restrictions to each part form
a triangle in $A_0$ and a matching in $A_1$.
There are $6$ lines of type $A_0^3 A_1$, and we note that
the $3$ points of $A_0$ consist of $1$ point from $\{x,y,z\}$
and $2$ of the other $4$ points.
This leaves $4$ lines of type $A_0 A_1^3$, for which we note
that the point in $A_0$ is not one of $\{x,y,z\}$.
In the picture, we have
$x=100$, $y=010$, $z=001$, $A_0 = $ white discs, $A_1 = $ black discs.
}
\parbox{0.45\textwidth}{
\includegraphics[height=7cm]{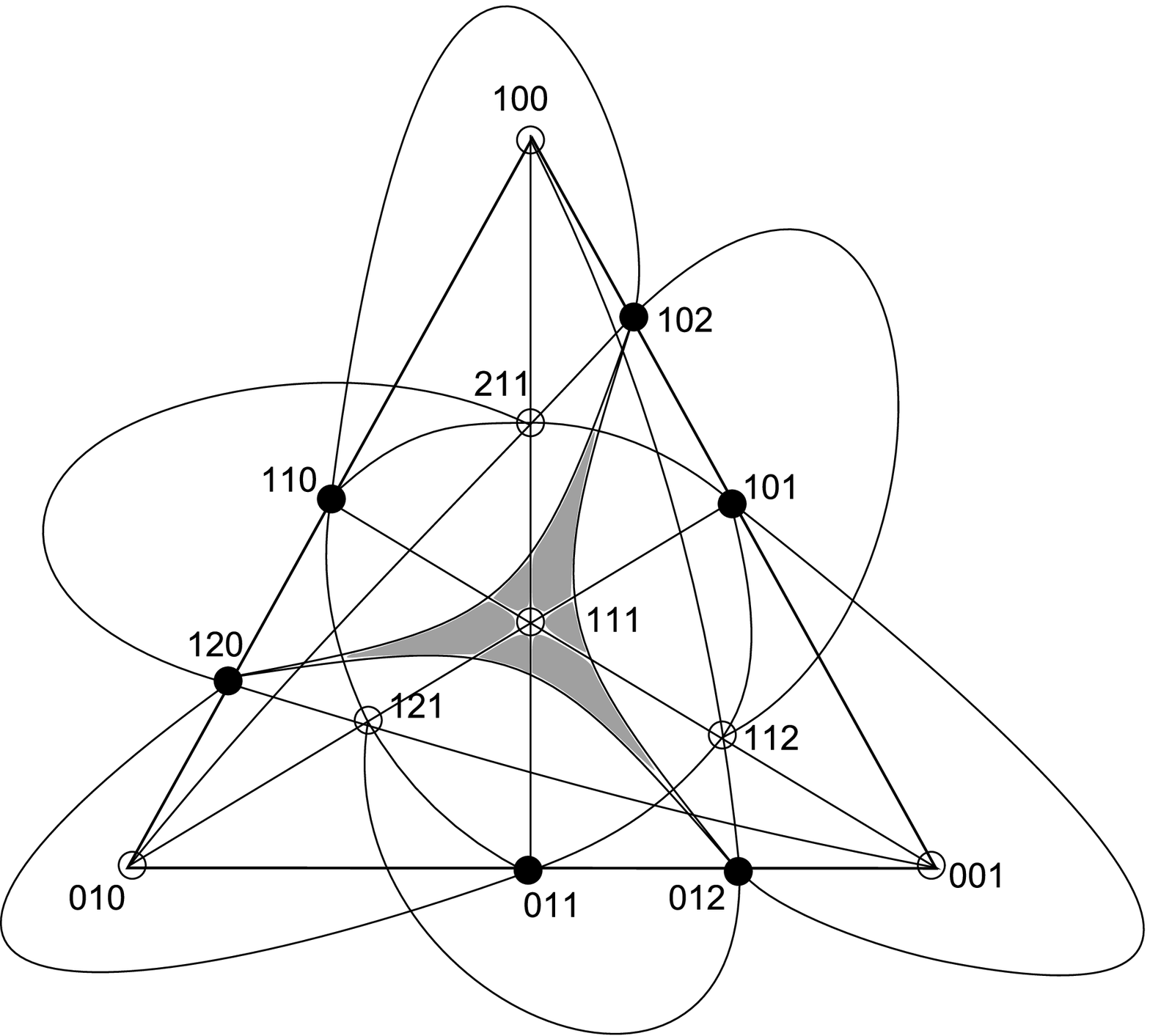}
}

\nib{Construction.} Now we can describe an improvement to the complete oddly bipartite construction.
Take a set $X$ of $n$ points, $n$ even, and partition it
as $X=X_0 \cup X_1$ with $|X_0|=n/2+1$. Choose $2$ special points
$a,b$ in $X_0$. Form a $4$-graph $H$ by
taking as edges all $4$-tuples which either have $3$ points in one $X_i$ and
$1$ point in $X_{1-i}$, or have $2$ points in $X_1$ and $2$ points in $X_0$,
exactly one of which is $a$ or $b$. Then $H$ does not contain $PG_2(3)$,
as the edges with $2$ points in each part do not contain the triangle
with matching neighbourhood configuration described in the previous paragraph.
Also, the minimum $3$-degree of $H$ is $n/2-1$.

\medskip

We will show that this construction is optimal. First we need the following lemma,
which is similar to (but does not follow from) a case of
a result of Diwan and Sobhan Babu \cite{DS}.
The proof is a simple but slightly tedious case analysis which we outline here
for the reader's convenience.

\begin{lemma} \label{tcm}
Suppose $G_1$, $G_2$ and $G_3$ are graphs on the same set $Y$ of at least $8$
vertices, so that
each has minimum degree at least $2$ and there is no `multicoloured' matching
$e_1$, $e_2$, $e_3$ with $e_i \in G_i$ for each $i$.
Then there are two points $a,b$ that meet each edge of each $G_i$.
\end{lemma}

\nib{Proof.}
Suppose that no two points meet each edge of each $G_i$. We claim that there
is a matching $M$ of size $3$ in $G = G_1 \cup G_2 \cup G_3$.
To see this, consider a possible counterexample $G$.
Divide into cases according to the connectivity $\kappa(G)$ of $G$.
The minimum degree condition implies that each component of $G$ contains a matching
of size $1$ (i.e. an edge), and even a matching of size $2$ unless it is
a triangle, so $G$ must be connected. If $\kappa(G)=1$ and $v$ is a cutvertex
then every component of $G \sm v$ contains an edge, so $G \sm v$ has exactly
two components $C$, $C'$. Since there is no matching of size $3$, at least one,
say $C$, is a star with at least $4$ vertices,
i.e. its edges all contain some vertex $x$ in $C$.
By the minimum degree condition any other vertex $y$ in $C$ is joined to $v$.
But now we can find a matching of size $3$: take $ya$, $xz$ for some $z \in C \sm \{x,y\}$
and an edge in $C'$. Next suppose that $\kappa(G)=2$ and $\{u,v\}$ is a cutset.
By hypothesis $\{u,v\}$ does not meet every edge, so some component $C$
of $G \sm \{u,v\}$ contains an edge. If there are at least two other components $C'$, $C''$
of $G \sm \{u,v\}$ then we can extend this to a matching of size $3$ with
an edge from $C$ to $u$ and an edge from $C''$ to $v$. If there is one other component $C'$
of $G \sm \{u,v\}$ then we can find a matching of size $3$ in which there is an edge
from $C$ to $u$ or $v$, from $C'$ to the other of $u$ or $v$, and a third edge in $C$.
Finally, suppose $\kappa(G) \ge 3$ and $S$ is a minimum cutset. Since every $x \in S$
has a neighbour in every component of $G \sm S$ there are exactly $2$ components of $G \sm S$.
Then we can find a matching of size $3$ using two edges from $S$ to components of $G \sm S$
and one edge inside a component of $G \sm S$. In all cases we see that there
is a matching $M$ of size $3$ in $G = G_1 \cup G_2 \cup G_3$.

If not all three colours $G_i$ are used in $M$ we will show how to increase
the number of colours. Suppose first that all three edges are in $G_1$.
Pick any new point $x$ and an edge $e$ of $G_2$ that contains it.
We can include $e$ and discard whichever
edge of our original matching meets it (or any if none does)
to obtain a new matching $M'$ of size $3$ on
which the colour $G_2$ also appears (on edge $e$). Now
any new point $y$ must be incident to exactly two edges of
$G_3$, joining it to the endpoints of $e$ (otherwise we would
have a multicoloured matching). However, we could take another
new point $z$, an edge $e'$ in $G_2$ containing $z$, and
include $e'$ in a new matching (discarding an appropriate edge).
Now we have either $2$ edges of $G_1$ and $1$ of $G_2$ (if $e$
was discarded) or $2$ edges of $G_2$ and $1$ of $G_1$ (otherwise).
Either way, the same reasoning as before tells us
that every new point $y'$ is incident in $G_3$ to exactly to the
endpoints of $f$, where $f \ne e$ is some edge of the new
matching, contradicting the fact that it is incident exactly
to the endpoints of $e$. \qed

\begin{theo} \label{pg23-exact}
If $n$ is sufficiently large and $H$ is a $4$-graph on a set $X$ of $n$ vertices
that does not contain $PG_2(3)$ then $\delta_3(H) < n/2$.
\end{theo}

\nib{Proof of Theorem \ref{pg23-exact}.}\
Suppose for a contradiction that $\delta_3(H) \ge n/2$.
By Theorem \ref{pg2-pre-exact} we have $n$ even, $X = X_0 \cup X_1$
with $|X_0|=|X_1|=n/2$ and $H[X_0]=H[X_1]=\emptyset$. Also, for every
$3$-tuple contained in one of the $X_i$, to get minimum $3$-degree $n/2$
every $4$-tuple obtained by adding a point from $X_{1-i}$ must be an edge.
For every pair $a,b$ in $X_0$ we have a graph $G_{a,b} = N_H(a,b) \cap \binom{X_1}{2}$
on $X_1$ with minimum degree at least $2$. These graphs do not contain
a `triangle-coloured matching', i.e. a triple $a,b,c$ in $X_0$ and
a matching $e_1, e_2, e_3$
in $X_1$ with $e_1 \in G_{b,c}$, $e_2 \in G_{a,c}$, $e_3 \in G_{a,b}$.
For every $4$-tuple with $3$ points in one $X_i$ and one point in $X_{1-i}$ is an edge, so
using the description of the blocking sets in $PG_2(3)$,
any triangle-coloured matching could be
completed to a copy of $PG_2(3)$, contrary to assumption.
By the lemma there must be two points $a_1,b_1$ in $X_1$ that meet
every edge of each of $G_{b,c}$, $G_{a,c}$ and $G_{a,b}$. In fact $a_1,b_1$
must meet every edge of $G_{a',b'}$ for any pair $a',b'$ in $X_0$,
as may be seen by applying the previous reasoning in the triangles
$a' b' a$ and $b' a b$ (without loss of generality). Similarly there
are two special points $a_0,b_0$ in $X_0$ that meet every edge
with $2$ points in each $X_i$. But now any triple
$a,b,c$ with say $a,b$ in $X_0 \sm \{a_0,b_0\}$ and
$c$ in $X_1 \sm \{a_1,b_1\}$ is not contained in any edge with
$2$ points in each $X_i$, so has $3$-degree equal to $n/2-2$.
This contradiction completes the proof. \qed

\subsection{The projective plane over $\mb{F}_4$}

We conclude by demonstrating a somewhat surprising phenomenon for
the projective plane over $\mb{F}_4$: its codegree density is
less than $1/2-c$ for some absolute $c>0$, unlike the cases
of $\mb{F}_2$ and $\mb{F}_q$, $q$ odd where the codegree density is $1/2$.

Before proving our bound we need some information about
the blocking sets of $PG_2(4)$. A classification was given in
\cite{BE}, but we will just need two specific examples.
For the first, suppose more generally that $q$ is a prime power.
Then one blocking set in $PG_2(q^2)$ is a Baer subplane $B$,
which may be constructed by restricting to those
points $(x_0:x_1:x_2)$ that have some representative
$(a_0,a_1,a_2)$ with each $a_i$ in the base field $\mb{F}_q$.
Each line of $PG_2(q^2)$ contains either
$1$ or $q+1$ points of $B$. The intersections of
size $q+1$ in $B$ form the lines of a copy of $PG_2(q)$. Also,
since every pair of lines in $PG_2(q)$ intersect, the
lines of $PG_2(q^2)$ containing them do not intersect
outside of $B$. For the sake of being more concrete, we remark
that this can be described more explicitly in $PG_2(4)$ using
the representation $PG_2(4) = \{ A + x : x \in \mb{Z}_{21}\}$, where
$A = \{3,6,7,12,14\}$. An example of a Baer subplane is
$B = \{x : x \equiv 0 \mbox{ mod } 3\}$. Dividing by $3$ we can
represent the lines of $B$ as $\{ A' + x: x \in \mb{Z}_7 \}$,
where $A'=\{1,2,4\}$: a well-known description of the Fano plane.

\begin{figure}
\begin{center}
\includegraphics[height=6cm]{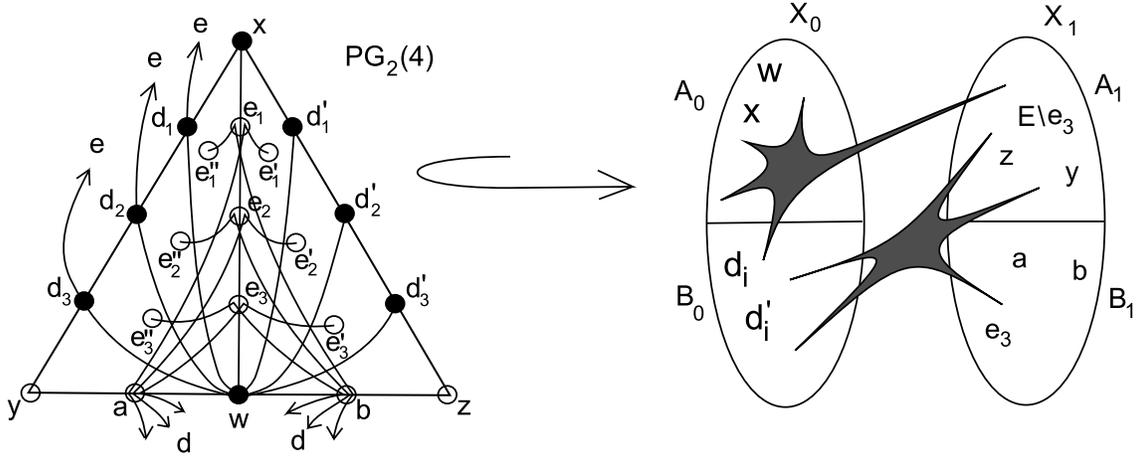}
\end{center}
\caption{$PG_2(4)$ and the second embedding used in the proof.}
\label{pg24-pic}
\end{figure}

We also need to consider the blocking set \d{$\wedge$}
$= L(x,y) \cup L(x,z) \cup \{w\} \sm \{y,z\}$,
where $x,y,z$ are $3$ non-collinear points, and $w \in L(y,z)$.
Consider the associated $2$-colouring
$C_0 =$ \d{$\wedge$}, $C_1 = V(PG_2(4)) \sm C_0$.
Recall that the type of an edge $L$ is $C_0^{|L \cap C_0|} C_1^{|L \cap C_1|}$.
There are $2$ edges of type $C_0^4 C_1^1$,
$9$ of type $C_0^1 C_1^4$,
$3$ of type $C_0^3 C_1^2$ and
$7$ of type $C_0^2 C_1^3$.
For the $C_0^3 C_1^2$ edges the
triples in $C_0$ can be described as $wd_id_i'$ with $i=1,2,3$, where
$L(x,y) = \{x,y,d_1,d_2,d_3\}$ and
$L(x,z) = \{x,z,d_1',d_2',d_3'\}$.
For the $C_0^2 C_1^3$ edges the triples in $C_1$ can be
described as  $e_1e_2e_3$, $ae_i e_i'$, $be_ie_i''$ with $i=1,2,3$,
where we
write $L(y,z)=\{y,z,w,a,b\}$, $L(x,w) = \{x,w,e_1,e_2,e_3\}$,
$L(a,e_i) \sm (L(x,y) \cup L(x,z) \cup \{a,e_i\}) = \{e_i'\}$ and
$L(b,e_i) \sm (L(x,y) \cup L(x,z) \cup \{b,e_i\}) = \{e_i''\}$.
Figure \ref{pg24-pic} shows an incomplete representation of $PG_2(4)$,
with $C_0$ represented by black discs and $C_1$ by white discs.

The following statement generalises Theorem
\ref{pg24}, in that it makes an allowance for a small exceptional set
of small codegrees.

\begin{theo} \label{gen-pg24}
For any $\eps>0$ sufficiently small there
is $\theta>0$ and $n_0$ so that if $H$ is a
$5$-graph on a set $X$ of $n>n_0$ vertices
and $G = \{Q \in \binom{X}{4}: |N_H(Q)| <  (1/2-\theta) n \}$
satisfies $d(G) < \theta$ then $d_{PG_2(4)}(H) > \eps$.
\end{theo}

\nib{Proof.}
Let $\eps_5$ be a small constant (say $10^{-6}$),
introduce a hierarchy of constants
$\eps_5 \gg \eps_4 \gg \cdots \eps_3 \gg \eps_2
\gg \eps_1 \gg \eps_0 \gg \eps \gg \theta$
and suppose $n$ is sufficiently large.
Suppose for a contradiction that
$d(G) < \theta$ but $d_{PG_2(4)}(H) < \eps$.
By Theorem \ref{pg2-stab} we can find
a partition $X = X_0 \cup X_1$
so that $|X_0|$ and $|X_1|$ are $(1/2 \pm \eps_0)n$
and $|H[X_0]|+|H[X_1]| < \eps_0 n^5$.

We introduce two $3$-graphs for $i=0,1$:
$$J_i = \left\{ e \in \binom{X_i}{3}: \left| N_H(e) \cap \binom{X_{1-i}}{2} \right|
>\eps_1\binom{|X_{1-i}|}{2} \right\}.$$
We claim that $d_{PG_2(2)}(J_i) < \eps_2$. For suppose otherwise, say
that $d_{PG_2(2)}(J_0) > \eps_2$, and let $\phi$ be a random map from
$V(PG_2(4))$ to $X$, conditioned on the event $E(\phi)$ that a Baer subplane $B$ is
mapped to $J_0$, and the other points are mapped to $X_1$.  We
estimate the probability that $\phi$ embeds $PG_2(4)$ in $H$.  There
are $2$ types of edges to consider: $14$ of type $X_0^1X_1^4$ and $7$
of type $X_0^3 X_1^2$. Suppose $L$ is an edge of $PG_2(4)$ that we are
attempting to embed with type $X_0^1 X_1^4$. For each point $a$ in
$X_0$ let $m(a)$ be the number of $4$-tuples $Q \in \binom{X_1}{4}$
such that $a \cup \{Q\}$ is not an edge of $H$. We have
\begin{eqnarray*}
(1/2-\theta)n \left(\binom{|X_1|}{4} - \theta \binom{n}{4} \right)
& < & \sum_{Q \in \binom{X_1}{4}} |N_H(Q)|
= \sum_{a \in X_0} \left| N_H(a) \cap \binom{X_1}{4} \right| + 5|H[X_1]| \\
& < & |X_0| \binom{|X_1|}{4} - \sum_{a \in X_0} m(a) + 5\eps_0 n^5,
\end{eqnarray*}
so
$\sum_{a \in X_0} m(a) < 10\eps_0 n^5$.  Now
$$\mb{P}(\phi(L)
\notin H) = \sum_{a \in X_0} \mb{P}(\phi(L) \cap X_0 = \{a\})
m(a)|X_1|^{-4},$$
and since there are at most $\binom{|X_0|-1}{6}$
copies of $PG_2(2)$ in $J_0$ that use $a$ we have $\mb{P}(\phi(L) \cap
X_0 = \{a\}) < \binom{|X_0|-1}{6} / \eps_2 \binom{|X_0|}{7} =
\frac{7}{\eps_2|X_0|}$, so $\mb{P}(\phi(L) \notin H) < \frac{7}{\eps_2|X_0|} \cdot
10\eps_0 n^5 |X_1|^{-4} < 5000\eps_2^{-1} \eps_0$.  We are
attempting to embed $14$ edges of type $X_0^1X_1^4$, so the
probability that any fails is less than $10^5 \eps_2^{-1} \eps_0$.
On the other hand, we have conditioned on the event that the Baer
subplane $B$ is mapped to $J_0$, and so by definition each of the $7$
attempted embeddings of type $X_0^3 X_1^2$ is successful with
probability at least $\eps_1 - O(1/n)$. Furthermore, if
$T_1,\cdots,T_7$ are the triples of the Baer subplane then
$N_{PG_2(4)}(T_1),\cdots,N_{PG_2(4)}(T_7)$ are disjoint sets, so the
events are independent. With probability $1-O(1/n)$ their images under
$\phi$ are disjoint, so we estimate $\mb{P}(\phi(PG_2(4)) \sub H) >
\eps_1^7 - O(1/n) - 10^5\eps_2^{-1} \eps_0 > \eps^{1/2}$, say. Then,
letting $\psi$ be a random map from $V(PG_2(4))$ to $X$,
we have $d_{PG_2(4)}(H) = \mb{P}(\psi(PG_2(4)) \sub H)
\ge \mb{P}[E(\psi)] \eps^{1/2}
> \eps_2(1/4)^{14}\eps^{1/2} > \eps$.
This contradiction shows that $d_{PG_2(2)}(J_i) < \eps_2$.

Next we introduce $2$-graphs for $i=0,1$:
$$P_i = \left\{ \{a_i,a_i'\} \in \binom{X_i}{2}
\mbox{ with } |N_{J_i}(a_i,a_i')| < (1/2-\eps_2)|X_i| \right\}.$$
We claim that $|P_i| < \eps_2 \binom{|X_i|}{2}$ for at least one of
$i=0,1$.  For suppose otherwise, and consider a $4$-tuple $Q =
\{a_0,a_0',a_1,a_1'\}$, where $\{a_i,a_i'\}$ are chosen uniformly at
random from $P_i$, $i=0,1$.  We have $\mb{E}|N_H(Q)| = \mb{E}|N_H(Q)
\cap X_0| + \mb{E}|N_H(Q) \cap X_1|$. Now
$$\mb{E} |N_H(Q) \cap X_0| =
\mb{E}_{a_0,a_0'} \sum_{a \in X_0} \mb{P}_{a_1,a_1'}(a \in N_H(Q) |
a_0,a_0')$$
and $\mb{P}_{a_1,a_1'}(a \in N_H(Q) | a_0,a_0') =
\mb{P}(\{a_1,a_1'\} \in N_H(a_0,a_0',a))$.  For $a \in
N_{J_0}(a_0,a_0')$ we estimate this probability trivially by $1$, but
for $a \in X_0 \sm N_{J_0}(a_0,a_0')$ we use the definition of $J_0$
and the lower bound on $P_1$ to estimate $\mb{P}(\{a_1,a_1'\} \in
N_H(a_0,a_0',a)) < \eps_1 \eps_2^{-1}$. This gives $\mb{E} |N_H(Q)
\cap X_0| \le |N_J(a_0,a_0')| + \eps_1 \eps_2^{-1} |X_0| <
(1/2-\eps_2/2)|X_0|$. Similarly we estimate $\mb{E} |N_H(Q) \cap X_1|
< (1/2-\eps_2/2)|X_1|$, so $\mb{E}|N_H(Q)| < (1/2-\eps_2/2)n$.
However, we also have
$\mb{E}|N_H(Q)| > \mb{P}(Q \notin G) \cdot (1/2-\theta)n$,
and $\mb{P}(Q \in G) \le \frac{|G|}{|P_0||P_1|} < \frac{\theta n^4/24}{
\eps_2\binom{|X_0|}{2}\eps_2\binom{|X_1|}{2}} < \eps$,
so $\mb{E}|N_H(Q)| > (1-\eps)(1/2-\theta)n$.  This
contradiction shows that at least one $P_i$ is small, say $|P_0| <
\eps_2 \binom{|X_0|}{2}$.

Now we can apply Theorem \ref{pg2-stab}
to find a partition $X_0 = A_0 \cup B_0$ where $|A_0|$ and $|B_0|$ are
$(1/2 \pm \eps_3)|X_0|$ such that at most $\eps_3 |X_0|^3$ edges of $J_0$
are contained entirely within $A_0$ or within $B_0$.
Next we repeat the argument to deduce similar structural information on $J_1$.
Let
$$P'_0 = \left\{ \{a_0,a_0'\} \in \binom{A_0}{2}
\mbox{ with } |N_{J_0}(a_0,a_0') \cap A_0| < 400\eps_3|X_0| \right\},$$
$$P'_1 = \left\{ \{a_1,a_1'\} \in \binom{X_1}{2}
\mbox{ with } |N_{J_1}(a_1,a_1')| < (1/2-\eps_4)|X_1| \right\}.$$
We must have $|P_0'| > \frac{1}{2}\binom{|A_0|}{2}$, otherwise
$$\eps_3 |X_0|^3  > |(J_0)_{A_0}|
= \frac{1}{3} \sum_{a_0,a_0' \in A_0} |N_{J_0}(a_0,a_0') \cap A_0|
> \frac{1}{6}\binom{|A_0|}{2} \cdot 400\eps_3|X_0|,$$
which is a contradiction.
Now we cannot have
$|P_1'| > \eps_4 \binom{|X_i|}{2}$, as then considering $Q =
\{a_0,a_0',a_1,a_1'\}$, where $\{a_i,a_i'\}$ are chosen uniformly at
random from $P_i'$, $i=0,1$, we estimate (similarly to before)
$\mb{E} |N_H(Q) \cap X_0| < (1/2+500\eps_3)|X_0| + \eps_1 \eps_4^{-1} |X_0|$
and
$\mb{E} |N_H(Q) \cap X_1| < (1/2-\eps_4)|X_1| + 10\eps_1 |X_1|$, so
$\mb{E} |N_H(Q)| < (1/2-\eps_4/2)n$, contradiction.  Again, by Theorem
\ref{pg2-stab} we find a partition $X_1 = A_1 \cup B_1$ where
$|A_1|$ and $|B_1|$ are $(1/2 \pm \eps_5)|X_1|$ such that at most
$\eps_5 |X_1|^3$ edges of $J_1$ are contained entirely within $A_1$ or
within $B_1$.

Let $T_i$ count edges of $H$ of type $A_0^2 A_1^2 B_i^1$, $i=0,1$.
We can bound $T_i$ by summing degrees of quadruples
$Q = \{a_0,a_0',a_1,a_1'\}$ with $\{a_i,a_i'\} \in \binom{A_i}{2}$:
$$(1/2-\theta)n \left(\binom{|A_0|}{2} \binom{|A_1|}{2} - \theta n^2/24\right)
< T_0 + T_1 + 3 \sum_{e \in \binom{A_0}{3}} \left|N_H(e) \cap \binom{A_1}{2}\right|
+ 3 \sum_{e \in \binom{A_1}{3}} \left|N_H(e) \cap \binom{A_0}{2}\right|.$$
Now
$$\sum_{e \in \binom{A_0}{3}} \left|N_H(e) \cap \binom{A_1}{2}\right|
< |(J_0)_{A_0}| \binom{|A_1|}{2} + \binom{|A_0|}{3} \cdot \eps_1 \binom{|X_1|}{2}
< 5\eps_3 n \binom{|A_0|}{2}\binom{|A_1|}{2}$$
and
$$\sum_{e \in \binom{A_1}{3}} \left|N_H(e) \cap \binom{A_0}{2}\right|
< |(J_1)_{A_1}| \binom{|A_0|}{2} + \binom{|A_1|}{3} \cdot \eps_1 \binom{|X_0|}{2}
<  5\eps_5 n \binom{|A_0|}{2}\binom{|A_1|}{2},$$
so
$T_0 + T_1 > (1/2-40\eps_5)n \binom{|A_0|}{2} \binom{|A_1|}{2}$.
Therefore
\begin{align*}
T_0 & > |B_0| \binom{|A_0|}{2} \binom{|A_1|}{2}
- \left( (|B_0|+|B_1|)  \binom{|A_0|}{2} \binom{|A_1|}{2}
- (1/2-40\eps_5)n \binom{|A_0|}{2} \binom{|A_1|}{2}
\right) \\
& > (1-400\eps_5) |B_0| \binom{|A_0|}{2} \binom{|A_1|}{2}.
\end{align*}
By symmetry, similar bounds hold for the number of edges in each
case when we specify a triple in one $X_i$ respecting
the partition $(A_i,B_i)$ and a pair in $A_{1-i}$ or $B_{1-i}$,
i.e. the types
$A_0^2 B_0^1 A_1^2$, $A_0^2 B_0^1 B_1^2$,
$A_0^1 B_0^2 A_1^2$, $A_0^1 B_0^2 B_1^2$,
$A_1^2 B_1^1 A_0^2$, $A_1^2 B_1^1 B_0^2$,
$A_1^1 B_1^2 A_0^2$, $A_1^1 B_1^2 B_0^2$.

Now we find $PG_2(4)$ using the  \d{$\wedge$}
colouring. Let $\phi$ be a random map
from $V(PG_2(4))$ to $X$, conditioned on the event $E'(\phi)$ that
$w,x$ are in $A_0$,
$y,z,e_1,e_2, e_i', e_i''$ are in $A_1$,
$d_i,d_i'$ are in $B_0$
and $e_3,a,b$ are in $B_1$
($i$ ranges from $1$ to $3$).
Note that $M = \min\{|A_0|,|A_1|,|B_0|,|B_1|\} > (1/4-2\eps_5)n$.
There are $11$ attempted embeddings of type $X_i^4 X_{1-i}^1$.
Recalling that $m(a)$ is the number of $4$-tuples $Q \in \binom{X_1}{4}$
such that $a \cup \{Q\}$ is not an edge of $H$,
we estimate that
the probability that we fail to embed some such $L$
is at most
$11 \sum_{a \in X_{1-i}} m(a) M^{-5} < 10^6 \eps_0$.
There are $10$ attempted embeddings
that have the type discussed in the previous paragraph,
i.e. one of the types equivalent to $A_0^2 B_0^1 A_1^2$.
Each fails with probability at most $400\eps_5$, so the probability
that any fails is at most $4000\eps_5$.
Now a random map $\psi$ from $V(PG_2(4))$ to $X$
satisfies $E'(\psi)$ with probability at least $(1/5)^{21}$
and so succeeds in embedding $PG_2(4)$ in $H$
with probability at least $(1/5)^{21}/2 > \eps$.
This contradiction completes the proof. \qed

\medskip

\nib{Remark.} The question of what happens for $PG_2(2^s)$ in
general is intriguing. It seems plausible that the above
approach of going from $PG_2(2)$ to $PG_2(4)$ could be adapted
to an inductive argument when $s=2^r$ is a power of $2$.
Much of the argument would go through as above: our hypergraph
has an approximate bipartition $X=X_0 \cup X_1$ and
$J_i = \{e \in \binom{X_i}{2^{2^{r-1}}+1}: |N_H(e)|
> \eps_1\binom{|X_{1-i}|}{2^{2^r}-2^{2^{r-1}}} \}$ satisfy
$d_{PG_2(2^{2^{r-1}})}(J_i) < \eps_2$. The step that may fail is
finding $Q$ of low degree: our approach used
the convenient coincidence of
$2 \cdot 2^{2^{r-1}} = 2^{2^r}$, which only occurs for $r=1$.

\section{Future directions}

The basic form of our hypergraph regularity method has been well illustrated
by its application to projective planes, which are relatively easy to deal
with (for reasons yet to be understood), although even here we cannot give
exact answers in all cases, and fields of even size seem particularly strange.
However,
the quasirandom counting lemma has potential to be a powerful tool
in the study of any Tur\'an problem, whether generalised or standard.
For example, in the Tur\'an problem
for the tetrahedron, if we consider any $K^3_4$-free $3$-graph $H$ and a
vertex $x$ then the edges of $H$ cannot be quasirandomly distributed with
positive density within the triangles of the neighbourhood graph $N_H(x)$.

If we restrict attention only to excluding simple
$k$-graphs $F$ (meaning that each pair of edges in $F$ have at most one
common point) then the projective geometries in higher dimensions point
to one stumbling block that should be overcome in future developments of this
method. For example, if we consider a $3$-graph $H$ on $n$ vertices
with no $PG_3(2)$ and
all but $o(n^2)$ codegrees at least $(2/3-o(1))n$ then our results
will give a set $Z$ of $(2/3-o(1))n$ vertices that induce a $3$-graph with
no Fano plane and all but $o(n^2)$ codegrees at least $(1/2-o(1))|Z|$,
so $Z$ is approximately bipartite by our structure result. This suggests
that $PG_3(2)$ should be approximately tripartite, and hence an inductive
approach for general $m$ and $q$ showing that a $(q+1)$-graph $H$
with no $PG_m(q)$ and all but $o(n^q)$ codegrees at least $(1-1/m-o(1))n$
should be approximately $m$-partite.

A potential approach to filling in
the gap is suggested by the remark after Theorem \ref{structure}. If we
stick to $PG_3(2)$ for the sake of simplicity, then not only do we have
an approximately bipartite subhypergraph of size about $(2/3)n$, but
any set of vertices $V'$ obtained by taking some classes of the regularity
partition contains some approximately bipartite $Z'$ of size about
$(2/3)|V'|$. Thus we are faced with the problem of recovering structural
information about $H$ from various restrictions, which is perhaps
best understood in the context of property testing (see \cite{RSc2} for
a hypergraph property testing result and many references to the literature).
Although a full investigation
of this idea is beyond the scope we have set for this paper introducing
our basic method, we remark that
it should be possible to carry our arguments over to this context via
a random reducibility property of quasirandom complexes, i.e. that a random
restriction of a quasirandom complex to sets of large constant size
should be quasirandom with high probability:
a high-level sketch is that the martingale
used by Lov\'asz and Szegedy \cite{LS} to show concentration
of the probability that a random map from a fixed graph $F$ to a random
graph $G$ is a homomorphism may be extended to show concentration
of the octahedral counting function that appears in the definition
of quasirandomness. This will allow
us to conclude that if $m$ is a large constant and $M$ is a random
$m$-set of vertices then $M$ contains an approximately bipartite
subhypergraph of size about $(2/3)m$ with high probability
(say $1 - \exp -m^c$ for some $c>0$). However, even assuming this it is
still not clear how to recover the global approximate
structure of $H$. If exact results are desired we also have the
problem of recovering the exact structure from the approximate structure.
This seems to be quite a different type of question, and so far all
instances of its solution have been of a rather ad hoc character, so
it would be interesting to develop some general principles here as well.

\medskip

\nib{Acknowledgement.} The author thanks Oleg Pikhurko and Yi Zhao
for discussions about Mubayi's conjecture and the anonymous referees
for helpful comments.

\appendix

\section{Variant forms of the Gowers quasirandomness framework}

In this appendix we justify the variant forms of the decomposition theorem
and counting lemma that we used in the paper. We first note that the proof
given for Theorem 5.1 in \cite{G2} also proves Theorem \ref{gencount}
in this paper: it only helps to assume stronger inequalities for the
$\eta_i$, and our parameter hierarchy is such that we can replace
$|\mc{J}_0 \sm \mc{J}_1| \eps$ in \cite{G2} by $\eps$ for simpler notation.
Then the arguments in section 6 of \cite{G2} go through as written.
The proof of Theorem \ref{decomp} is quite similar to the proof of Theorem 7.3
in \cite{G2}. We will not reproduce that proof in full, but will
outline it to sufficient extent to explain what modifications
are needed. First we recall a definition (\cite{G2} p. 36):

If $S = \{S_1,\cdots,S_s\}$ and $T = \{T_1,\cdots,T_t\}$ are partitions
of the same set $U$, the mean-square density of $S$ with respect to $T$
is
$$\mbox{msd}_T(S) = \sum_{i=1}^s \sum_{j=1}^t \frac{|T_j|}{|U|} \left(
\frac{|S_i \cap T_j|}{|T_j|} \right)^2.$$
We also recall Lemma 8.1 of \cite{G2}, which states that if $T'$ is
a refinement of $T$ then $\mbox{msd}_{T'}(S) \ge \mbox{msd}_T(S)$.

\medskip

\nib{Proof of Theorem \ref{decomp}.}
The proof of Theorem 7.3 in \cite{G2} is by means of the following iterative
procedure. Suppose we have a partition $k$-system $P$ (which
may be the initial partition or one produced by a some number of
iterations: we will not complicate notation with a sequence
$P_1,P_2,\cdots$). Consider a partition $k$-system $P^*$ defined by
what is called weak equivalence in \cite{G2}:
$S, S' \in K_A(X)$ are in the same class of $P^*_A$
exactly when $S_B$ and $S'_B$ are in the same class of $P_B$
for every proper subset $B \subset A$.
It is shown (see \cite{G2} pp. 38--39) that
if $\mb{P}_x(P(x)$ is $(\eps,j,k)$-quasirandom$)<1-\eps$,
then there is $A \in \binom{[r]}{\le k}$ and a refining
partition $k$-system $Q$ of $P$, so that
(a) $Q_B = P_B$ unless $B \subset A$, $|B|=|A|-1$,
(b) if $B \subset A$, $|B|=|A|-1$ then
$Q_B$ is a refinement of $P_B$ where each
class of $P_B$ is partitioned into at most $c_A(\eps,P)$
further classes, and 
(c) $\mbox{msd}_{Q^*_A}(Q_A) \ge \mbox{msd}_{P^*_A}(P_A)
+ f_A(\eps,P)$. Here $c_A$ and $f_A$ are explicitly defined functions
that depend only on $\eps$ and $P$, and furthermore
the dependence of $c_A$ on $P$ depends only $\{|P_B|: B \subset A\}$ and
that of $f_A$ depends only on $\{|P_B|: |B| \ge |A|\}$.
This property of $f_A$ implies that $f_A(\eps,Q)=f_A(\eps,P)$.
Furthermore, this argument still applies in the context of our proof, i.e.
we have the same conclusion
if $P(x)$ is $\ov{\eta}$-quasirandom and $\ov{d}$-dense
with probability less than $1-\eps$.
(Note that the functions $c_A$ and $f_A$ now depend on the functions
implicit in the $\ll$-notation for the parameter hierarchy.)

To see that the procedure terminates (with some system of
partitions with the required property) introduce a function
$\zeta_P$ for the system of partitions $P$, which is defined on
$\binom{[r]}{\le k}$ by
$\zeta_P(A) = \left\lceil \frac{1-\mbox{msd}_{P^*_A}(P_A)}{f_A(\eps,P)}
\right\rceil$.
Choose an ordering $<$ of $\binom{[r]}{\le k}$
in which $|B| \ge |B'|$ implies that $B<B'$. Order functions
$\zeta$ on $\binom{[r]}{\le k}$ by $\zeta<\zeta'$ if there is
$B \in \binom{[r]}{\le k}$ such that
$\zeta(B')=\zeta(B)$ for all $B'<B$ and $\zeta(B)<\zeta'(B)$.
This is a well-ordering, and the iteration takes the system $P$
to a system $Q$ with $\zeta_Q < \zeta_P$, so the procedure
terminates.

To prove our version we introduce further refinements
in each step of the procedure. First of all we make
the general observation that given any partition
$E = E_1 \cup \cdots \cup E_t$ there is an `equalising method'
to find a partition
$E = E'_0 \cup E'_1 \cup \cdots \cup E'_s$, for which
$|E'_i| = \lfloor |E|/t^2 \rfloor$, $|E'_0| < |E|/t$
and every $E'_i$ with $i \ne 0$ is contained in
some $E_j$. The method is to repeatedly and arbitrarily select
classes $E'_i$ within some $E_j$ that still has size at
least $\lfloor |E|/t^2 \rfloor$, and then remove
its elements from consideration in later stages. Thus
we are unable to use at most $|E|/t^2$ elements from each
of the $t$ original classes $E_i$, and we put the unused elements
together in an exceptional class
$E'_0 = E'_{0,1} \cup \cdots \cup E'_{0,t}$ of size at most $|E|/t$.

We start by using the equalising method in an
initial refinement to transform $P$ into some $Q$
with an equitable partition of the vertex set.
By first arbitarily refining $P$ we can assume that
$|X_i^0| < \eps|X_i|/2$, $1 \le i \le r$.
Then we repeatedly apply the same refinement procedure as above
followed by the equalising method:
at each stage we obtain a new system $Q$ from the original procedure
and then refine it to some equitable $Q'$.
It is clear that the number of classes remains bounded
by a function only of $m$, $r$, $k$ and $\eps$.
By Lemma 8.1 of \cite{G2} we have
$\mbox{msd}_{Q^{'*}_A}(Q'_A) \ge \mbox{msd}_{Q^*_A}(Q_A) \ge
\mbox{msd}_{P^*_A}(P_A) + f_A(\eps,P)$. Also, we still have
$Q_B=P_B$ for all $B$ with $|B| \ge |A|$, so
$f_A(\eps,Q)=f_A(\eps,P)$ and the iterations terminate as before.
The amounts added to the exceptional classes
decrease rapidly with each iteration (certainly
each is at most half of that at the previous
iteration), and as we initially added at most
$\frac{1}{2} \eps |X_i|$ exceptional elements of $X_i$
we end up with at most $\eps |X_i|$. \qed

\end{document}